\documentclass[a4paper]{article}
\usepackage[latin1]{inputenc}
\usepackage{amsthm}
\usepackage{amsfonts}
\usepackage[dvips]{graphicx}
\usepackage{latexsym}
\usepackage{amsmath}
\usepackage{color}

\usepackage{amssymb}

\newtheorem{thm}{Theorem}
\newtheorem{defn}{Definition}

\newtheorem{oss}{Remark}

\usepackage{color}
\definecolor{lightgray}{gray}{0.75}

\newcommand\greybox[1]{%
  \vskip\baselineskip%
  \par\noindent\colorbox{lightgray}{%
    \begin{minipage}{\textwidth}#1\end{minipage}%
  }%
  \vskip\baselineskip%
}
\newcommand{\tn}[1]{\mathbb{#1}}

\newcommand{\Krejci}{Krej\v c\'\i}
\DeclareMathOperator{\dive}{div}
\newcommand{\teta}{\theta}

\newcommand{\dt}{\partial_t}
\newcommand{\vc}[1]{{\bf #1}}
\newcommand{\vt}{\theta}
\newcommand{\bu}{{\bf u }}
\newcommand{\ub}{{\bf u }}
\newcommand{\vu}{{\bf u}}
\newcommand{\Grad}{\nabla_x}
\DeclareMathOperator{\deriv}{d}

\newcommand{\ddt}{\frac{\deriv\!{}}{\dit}}

\newcommand{\RR}{\mathbb{R}}

\allowdisplaybreaks[4]

\def\xibf{{\boldsymbol \xi}}
\def\phibf{{\boldsymbol \phi}}

\definecolor{green}{rgb}{0.0,0.4,0.2}
\definecolor{violet}{rgb}{0.4,0,0.9}
\numberwithin{equation}{section}

\newenvironment{bettirev}{\color{blue}}{\color{black}}
\newcommand{\bber}{\begin{bettirev}}
\newcommand{\eber}{\end{bettirev}}

\newenvironment{michelarev}{\color{red}}{\color{black}}
\newcommand{\bmicr}{\begin{michelarev}}
\newcommand{\emicr}{\end{michelarev}}

\begin{document}


\newcommand{\fhi}{\varphi}

\newcommand{\Vn}{{\bf V}_{\bf n}}
\newcommand{\rhs}{right hand side}
\newcommand{\lhs}{left hand side}
\newcommand{\io}{\int_\Omega}
\newcommand{\ito}{\int_0^t\int_\Omega}
\newcommand{\iTo}{\int_0^T\int_\Omega}
\newcommand{\iga}{\int_\Gamma}

\newcommand{\dix}{{\rm d} x}
\newcommand{\dis}{{\rm d} s}
\newcommand{\dit}{{\rm d} t}
\newcommand{\diy}{{\rm d} y}
\newcommand{\diS}{{\rm d} S}
\newcommand{\dir}{{\rm d} r}

\newcommand{\dn}{\partial_{\bf n}}

\newcommand{\baru}{\overline{u}}
\newcommand{\barf}{\overline{f}}
\newcommand{\barv}{\overline{v}}
\newcommand{\barz}{\overline{z}}
\newcommand{\barteta}{\overline{\theta}}
\newcommand{\bartheta}{\overline{\theta}}


\title{On a non-isothermal diffuse interface model for two-phase flows of incompressible fluids}

\author{
{Michela Eleuteri\thanks{Dipartimento di Matematica ``F. Enriques'',
Universit\`{a} degli Studi di Milano, Milano I-20133, Italy.
E-mail: \textit{Michela.Eleuteri@unimi.it}.  The author
is supported by the FP7-IDEAS-ERC-StG
Grant \#256872 (EntroPhase)}}\and
\and
{Elisabetta Rocca\thanks{Weierstrass Institute for Applied
Analysis and Stochastics, Mohrenstr.~39, D-10117 Berlin,
Germany. E-mail: \textit{rocca@wias-berlin.de} and Dipartimento di Matematica ``F. Enriques'',
Universit\`{a} degli Studi di Milano, Milano I-20133, Italy.
E-mail: \textit{elisabetta.rocca@unimi.it}.  The author
is supported by the FP7-IDEAS-ERC-StG
Grant \#256872 (EntroPhase)}}
\and
{Giulio Schimperna\thanks{Dipartimento di Matematica ``F. Casorati'', Universit\`a degli Studi di Pavia,
via Ferrata 1, Pavia I-27100, Italy.
E-mail: \textit{giusch04@unipv.it}. The author
is partially supported by the FP7-IDEAS-ERC-StG
Grant \#256872 (EntroPhase)}}}

\maketitle

\begin{abstract}\noindent
We introduce a diffuse interface model describing the evolution of a mixture of two
different viscous incompressible fluids of equal density. The main novelty of the present contribution consists in the fact that the effects of temperature on
the flow are taken into account. In the mathematical model,
the evolution of the velocity $\vu$ is ruled by the Navier-Stokes
system with temperature-dependent viscosity, while the order parameter $\fhi$ representing
the concentration of one of the components of the fluid is assumed to satisfy a
convective Cahn-Hilliard equation. The effects of the temperature are prescribed by a
suitable form of the heat equation. However, due to quadratic forcing terms, this equation
is replaced, in the weak formulation, by an equality representing energy
conservation complemented with a differential inequality describing production of entropy.
The main advantage of introducing this notion of solution 
is that, while the thermodynamical consistency is preserved, at the same time the energy-entropy formulation
is more tractable mathematically. Indeed, global-in-time existence for the initial-boundary value problem
associated to the weak formulation of the model is proved by deriving suitable a-priori
estimates and showing weak sequential stability of families of approximating solutions.
\end{abstract}

\smallskip


\smallskip


\section{Introduction}

We study a non-isothermal diffuse interface model for the flow of a mixture of two viscous incompressible
Newtonian fluids of equal density in a bounded domain $\Omega \subset \mathbb{R}^3$.
In classical models the interface between the two fluids is assumed to be a $2-$dimensional
sufficiently smooth surface; in this case capillarity phenomena are related to contact angle
conditions and a jump condition for the stress tensor across the interface. This classical
description fails when some parts of the interface merge or reconnect (developing singularities)
due to droplet formation or coalescence of several droplets. Indeed, despite the large amount of
mathematical literature on free boundary problems related to fluids with a classical
sharp interface, most papers are confined to the case of flows without singularities
in the interface and so far there is no satisfactory existence theory of weak solutions
for a two-phase flow of two viscous, incompressible, immiscible fluids with a classical sharp interface.
Thus, in order to avoid analytical problems related to interface singularities,
an alternative approach, based on {\it diffuse interface models}, can be used.
In this setting, the classical sharp interface, represented by a lower-dimensional
surface, is replaced by a thin interfacial region, whose ``thickness''
is described by a small parameter $\varepsilon > 0$. Therefore a partial mixing of
the macroscopically immiscible fluids is allowed; in order to describe this
phenomenon, a new variable $\varphi$ is introduced.
This quantity may represent the concentration difference or
the concentration of one component of the fluid. 

The original idea of diffuse interface model for fluids goes back to {\sc Hohenberg}
and {\sc Halpe\-rin} \cite{HH77} and it is usually referred with the name ``H-model''. Later,
{\sc Gurtin} et al.~\cite{GPV96} gave a continuum mechanical derivation based
on the concept of microforces. For a review of the development of diffuse interface
models and their applications we refer to \cite{Athesis07,AMW98} and the references
therein.

The present contribution is aimed at extending the {\sc H-model} to a non-isothermal setting,
developing a thermodynamically consistent theory. The system of partial differential
equations resulting from this approach couples the incompressible Navier-Stokes system
for the velocity ${\bf u}$ with a Cahn-Hilliard system  with convection, where we also
account for the effects of the (absolute) temperature $\theta$. Namely, we consider the
following equations:
\begin{align}
& \textnormal{div} \, {\bf u} = 0, \label{model1}\\[2mm]
& {\bf u}_t + {\bf u} \cdot \nabla_x {\bf u} + \nabla_x p = \dive \mathbb{S} - \varepsilon \textnormal{div}(\nabla_x \varphi \otimes\nabla_x \varphi), \qquad \mathbb{S} = \nu(\theta) D {\bf u} ,\label{model2}\\[2mm]
& \varphi_t + {\bf u} \cdot \nabla_x \varphi = \Delta \mu, \label{model3}\\[2mm]
& \displaystyle \mu = - \varepsilon \Delta \varphi + \frac{1}{\varepsilon} F'(\varphi)  - \theta, \label{model4}\\[2mm]
& c_V(\theta) \theta_t  + c_V(\teta) {\bf u} \cdot \nabla_x \theta+{\theta (\varphi_t+{\bf u}\nabla\varphi)} - \dive (\kappa(\theta) \nabla_x \theta) = \nu(\theta) |D {\bf u}|^2 + |\nabla_x \mu|^2 ,\label{model5}
\end{align}
in $\Omega \times (0, T)$, being $\Omega$ a bounded and sufficiently regular subset
of $\mathbb{R}^3$ and $T>0$ a given final time, which may be arbitrarily large.
Here $p$ is the pressure, $\mathbb{S} = \nu(\theta) D {\bf u}$ represents the dissipative
part of the stress tensor, where $\displaystyle D {\bf u} = (\nabla_x {\bf u} + \nabla^t_x {\bf u})/2$,
and $\nu(\theta)>0$ is the viscosity of the mixture. Moreover, $c_V(\theta)$ stands for the specific heat, $\kappa(\theta)$ indicates the heat conductivity,
$\varepsilon > 0$ is a (small) parameter related to the ``thickness'' of the interfacial region,
and $\mu$ is an auxiliary variable (usually named chemical potential) which helps particularly for the statement
of the weak formulation of the model. Finally, $F(\varphi)$ is some suitable energy density
whose expression is specified below. The capillarity forces due to surface tension are
modeled by an extra-contribution
$\varepsilon \nabla_x \varphi \otimes\nabla_x \varphi$
in the global stress tensor appearing in the right-hand side of \eqref{model1}.
System \eqref{model1}--\eqref{model5}
will be closed by adding the initial conditions and suitable boundary conditions.
Namely, the Cahn-Hilliard and temperature equations will be complemented by no-flux conditions,
while the velocity $\vu$ will be assumed to satisfy the so-called {\it complete slip}
conditions. As we will see (cf.~Subsection \ref{IBC} for more details), these choices
are crucial as we formulate the weak version of the model. Nevertheless,
some other choices could be considered as well (for instance, the case of periodic
boundary conditions can be treated similarly). It is also worth
noting that, integrating \eqref{model3} in space and using the no-flux
condition  together with \eqref{model1} and the complete slip
boundary condition for $\ub$, one gets back the mass conservation 
property (cf.~\eqref{zeromat} below).

Isothermal versions of our model have been studied by several authors
(see, e.g., \cite{A09,B,GG1,S,ZWH} and
references therein) and a mathematical theory can be now considered
to be well-established. On the
other hand, at least up to our knowledge, a non-isothermal model for two-phase fluids
has been analyzed only in the reference \cite{SLX09}, where a linearization of the internal
energy balance is used in order to describe the evolution of the temperature. This
permits the authors to get rid of the quadratic terms in
the \rhs\ of \eqref{model5}  and of the coupling beween \eqref{model4} and \eqref{model5}. 
On the other hand, the resulting model turns out to be thermodynamically consistent only 
in a neighbourhood of the equilibrium temperature. 

In this contribution, we describe the non-isothermal evolution of the fluid by
means of a model that keeps its thermodynamical consistency in a wide temperature range.
Moreover, we can prove global-in-time existence for a suitable weak formulation of
the associated initial-boundary value problem in three dimensions of space and
without any magnitude restriction on the data.
The analysis of non-isothermal problems in mathematical modelling of advanced materials
is gaining more and more importance in the recent years. We refer for instance,
without aiming at completeness, to \cite{Rou1,EKK12, EKK13, EKK14,Ulisse, FFRS, FRSZ, FR,  KR, KRS10, KRS09, KRS07, rocca-rossi1, rocca-rossi2,  RoRo, RR, RT, Rou2, Rou3},
where temperature-dependent  models are presented for describing
the evolution of several types of substances, like elastic media, plastic
materials (possibly with hysteresis, fatigue and damage), shape memory alloys, 
water-ice mixtures, and liquid crystals.
The idea of replacing the heat equation with the energy and entropy balances in
the weak formulation has been originally developed in~\cite{BFM09, FEI09}
in the framework of heat conduction phenomena in fluids and in \cite{FPR09}
in the case of solid-liquid phase transitions. It is worth observing that
the related notion of weak solution, which will be introduced in full detail
in Subsection~\ref{sec:pdes}, is consistent with the standard (strong) one.
Actually, it is not difficult to prove that, at least for {\it sufficiently smooth}\/
weak solutions, the total energy balance together with the entropy
inequality imply the original form of the heat (or, more precisely, internal energy
balance) equation \eqref{model5}. On the other hand, since this regularity
in our case is not at all known (for instance due to the occurrence of the 3D
Navier-Stokes system), this notion of solution turns out to be particularly useful 
because it allows us to prove a global in time existence result in 3D
and at the same time it guarantees the thermodynamical consistency of the model.
Better regularity properties are expected to hold for  weak solutions in the
2D case. This will be the subject of a forthcoming paper, where we will
also analyze the long-time dynamics of the model.

\paragraph{Plan of the paper.}  In the next Section~\ref{sec:model} we
provide a physical derivation of the model by following a variant of
the general approach devised by {\sc Fr\'emond} in \cite{Frem}.
Namely, the equations of the system are obtained by imposing the balances of energy
and entropy in terms of the free energy functional $\Psi$ and of the pseudopotential
of dissipation $\Phi$ and assuming standard consitutive relations. In the subsequent
Section~\ref{mainres} we introduce the main assumptions on data, which permit us to
state the weak formulation of the problem and the main existence theorem.
The proof of this result occupies the remainder of the paper and is split
into two steps: a-priori estimates, which are described in Section~\ref{sec:bou},
and weak-sequential stability, which is proved in the last Section~\ref{WSS}.


\section{Derivation of the model}
\label{sec:model}

We suppose that a two-component fluid occupies a bounded spatial domain
$\Omega\subset \RR^3$, with a sufficiently regular boundary $\Gamma$.
We let ${\bf n}$ denote the outer normal unit vector to $\Gamma$.
We denote by ${\bf u }= {\bf u}(t,x)$ the associated {\it velocity field}\/
in the Eulerian reference system.
Moreover, we introduce the {\em absolute temperature} $\theta(t,x)$
and the {\em order parameter} $\varphi(t,x)$, representing the concentration difference,
or the concentration of one component, of the fluid.
Furthermore, we denote as
$$
\frac{D w}{D t} = \dot{w} =w_t+{\bf u}\cdot\Grad w,
$$
the {\em material derivative} of a generic function $w$, while
$w_t$ (or also $\dt w$) stands for the partial derivative with respect to $t$.

We set $H := L^2(\Omega)$ and $V := H^1(\Omega)$. We will often write $H$ in place
of $H^3$, or $V$ in place of $V^3$, when vector-valued functions are
considered. In particular, $(\cdot,\cdot)$ will stand for the usual
scalar product both in $H$ and in $H^3$. For every $f \in V'$ we
indicate by $\overline{f}$ the {\it spatial mean} of $f$ over $\Omega$, i.e.
\[
  \overline{f} := \frac{1}{|\Omega|} \langle f, 1 \rangle,
\]
where $\langle \cdot, \cdot \rangle$ denotes the duality pairing between
$V'$ and $V$ and $|\Omega|$ stands for the Lebesgue measure of $\Omega$.
We note as $H_0$, $V_0$ and $V_0'$ the closed subspaces of functions
(or functionals) having zero mean value in $H$, $V$, and, respectively, in $V'$.
Then, by the Poincar\'e-Wirtinger inequality,
\[
  \| v \|_{V_0} := \bigg( \io | \Grad v |^2 \,\dix \bigg )^{1/2}
\]
represents a norm on $V_0$ which is equivalent to the norm inherited from $V$.
In particular $\|\cdot\|_{V_0}$ is a Hilbert norm and we can introduce the associated
Riesz isomorphism mapping $J:V_0 \to V_0'$ by setting, for $u,v\in V_0$,
\begin{equation}\label{defiJ}
  \langle J u, v\rangle:= (\!(u,v)\!)_{V_0}
   :=\int_\Omega\Grad u\cdot \Grad v\, {\rm d} x.
\end{equation}
For $f\in H_0$ it is easy to check that $u = J^{-1} f \in H^2(\Omega)$.
Actually, $u$ is the (unique) solution to the elliptic problem
\[
  u \in H_0, \quad
   -\Delta u = f, \quad
   \Grad u_{|_\Gamma} = 0.
\]
Moreover, if $u$ is as above, then
\[
  \big\langle J (u - \baru) , v \big \rangle
   = - \int_\Omega v \Delta u\, {\rm d} x
\]
for all $v \in V_0$. Finally, we can identify $H_0$ with $H_0'$
by means of the scalar product of $H$ so to obtain the Hilbert
triplet $V_0 \subset H_0 \subset V_0'$, where inclusions are continuous
and dense. In particular, if $z\in V$ and $v \in V_0$,
it is easy to see that
\begin{equation}\label{duale}
  \io \Grad z \cdot \Grad (J^{-1} v) \, \dix
   = \io ( z - \barz ) v \, \dix
   = \io z v \, \dix.
\end{equation}
%


\subsection{Free-energy and pseudopotential of dissipation}

\label{free-energy-sec}

We would like to apply here the general approach proposed in
the monograph \cite[Chapters 2, 3]{Frem} in order to
build our diffuse interface model for
incompressible fluids with conserved order parameter.

We start by introducing, in agreement with the basic
principles of classical Thermodynamics, the free-energy and the
pseudopotential of dissipation. To this aim, we specify
the set of the {\em state variables}, describing the
actual configuration of the material:
\begin{equation*}
  E=(\varphi,\Grad \varphi,\theta).
\end{equation*}
Correspondingly, the set of the {\em dissipative variables},
whose evolution describes the way along which
the system tends to dissipate energy, is given by
\begin{equation*}
  \delta E=\left(D{\bf u},\frac{D \varphi}{Dt}, \Grad\theta\right),
\end{equation*}
where
$$
   D{\bf u}:=\frac{\nabla_x {\bf u} + \nabla^t_x {\bf u}}{2}
$$
denotes the symmetric gradient of $\ub$.

Motivated by the Ginzburg-Landau theory for phase transitions, we choose the
free energy density $\psi$ and the free energy functional $\Psi$
in the form
\begin{equation}\label{psi}
  \psi(E) = \left(\frac{\varepsilon}{2}|\Grad\varphi|^{2}
  + \frac1\varepsilon F(\varphi)
  + f(\theta) - \teta\varphi\right), \qquad
  \Psi(E) = \io \psi(E) \, \dix,
\end{equation}
where $\varepsilon$ is a positive constant related to the interface thickness.
The function $F$ in \eqref{psi} penalizes the
deviation of the length $|\varphi|$ from its natural value 1; generally, $F$ is assumed to be
a sum of a dominating convex (and possibly non smooth) part
and a smooth non-convex perturbation of controlled growth.
Typical example are the standard {\em double-well} potential
$F(\varphi)=(|\varphi|^2-1)^2$ and the so-called {\em logarithmic
potential}\/
$F(\varphi)= (1+\fhi)\log (1+\fhi) + (1-\fhi)\log(1-\fhi) - \lambda \fhi^2$,
$\lambda\ge 0$.  Notice that in our case we will assume $F$ to be a double-well
potential and this hypothesis turns out to be essential in 
our analysis (cf.~Sec.~\ref{sec:bou}).  
Moreover, $f$ represents the purely caloric part of the free energy and
is linked to the specific heat $c_V(\teta)=Q'(\teta)$ by the
relation $Q(\teta)=f(\teta)-\teta f'(\teta)$ (cf.~\eqref{e} below). 
In what follows we will assume $c_V(\theta)\sim \theta^\delta$ with $\delta\in (1/2, 1)$. 
While the condition $\delta<1$ is quite natural (in particular it ensures the concavity 
of the entropy as function of $\theta$), the bound from below 
$\delta>1/2$ is crucial for the purpose of analyzing our PDE system 
(cf.~Sec.~\ref{sec:bou} below).

The evolution of the system is characterized by a second functional $\Phi$,
called {\em pseudopotential of dissipation}, assumed to be nonnegative
and convex with respect to the dissipative variables.
In order to present the explicit expression of $\Phi$, however,
we have to impose in some way the mass conservation constraint.
To this aim, we first decompose $\fhi = \fhi^0 + m_0$, where
$m_0$ represents the mean value of the initial datum $\fhi_0$.
Then, the conservation of mass corresponds to prescribe
$\fhi^0$ to take its values in~$H_0$ during the whole
evolution of the system.

Let us also note as $\Vn$ the subspace of $H^1(\Omega;\RR^3)$
consisting of the functions $\vu$ such that $\vu \cdot {\bf n}=0$
on $\Gamma$. Then, given a time-dependent family of divergence-free
vector fields $\vu(t,\cdot) \in \Vn$ and a scalar field
$\fhi=\fhi(t,x)$ satisfying the mass conservation constraint 
$\varphi(t,\cdot)=\varphi(0,\cdot)$ for a.e.~$t\in (0,T)$, 
we have
\begin{equation}\label{zeromat}
  \io \frac{D\fhi}{Dt}
   = \io ( \fhi_t + \vu \cdot \Grad \fhi ) \, \dix
   = 0.
\end{equation}
In other words, $\frac{D\fhi}{Dt}$ has zero spatial mean. Therefore, 
if $\fhi$ is so smooth to satisfy $\frac{D\fhi}{Dt}\in V_0'$
a.e.~in time, we can define
\begin{equation}\label{muzero}
  \mu^0 := - J^{-1} \frac{D\fhi}{Dt}, \qquad
   \text{so that } \frac{D\fhi}{Dt} = - J \mu^0 = \Delta \mu^0  \qquad\textnormal{in $V_0'$},
\end{equation}
and, consequently, $\mu^0\in V_0$. Hence, we can set
\begin{equation}\label{Phi}
  \Phi(\delta E,E)
   = \io \phi(\delta E, E) \, \dix
   + \left\langle \frac{D\fhi}{Dt}, J^{-1}\left(\frac{D\fhi}{Dt}\right) \right\rangle,
\end{equation}
where the ``local component'' $\phi$ of the
``dissipation density'' is given by
\begin{equation}\label{phi}
  \phi(\delta E,E) = \frac{\nu (\theta )}{2}|D{\bf u}|^{2}
   + I_{0}(\dive \ub)
   + \frac{\kappa(\theta )}{2\theta }|\Grad\theta |^{2}.
\end{equation}
Here, $\langle \cdot, \cdot \rangle$ denotes the duality
between  $V_0'$ and $V_0$,
$\nu = \nu(\theta) > 0$ is the viscosity coefficient
and $\kappa = \kappa(\theta) > 0$ represents the heat conductivity.
The {\em incompressibility} of the fluid is formally enforced by $I_0$, i.e.,
the indicator function of $\{0\}$ (given by $I_0 = 0$ if $\dive \ub =0$
and $+\infty$ otherwise).
The last term in \eqref{Phi}, which accounts for mass conservation,
is nonstandard and requires some words
of explanation. Basically, it corresponds to a (squared)
$V_0'$-norm of the material derivative
of $\fhi$ and, hence, depends on the dissipative
variable~$\frac{D\fhi}{Dt}$ in a {\it nonlocal} way. As will
be seen below, this gives rise to
a convective Cahn-Hilliard dynamics, in contrast
with the Allen-Cahn dynamics that would result
from the choice of the $H$-norm
(cf.~\cite{Frem} for more details).

The functionals $\Phi$ and $\Psi$ are assumed to be defined for all sets of
variables $E$ and $\delta E$ for which they make sense. In other words,
the finiteness, say, of $\Psi$ determines the class of admissible state variables.
So, in particular, if a time-dependent set of variables is given such that,
a.e.~in $(0,T)$, $\Psi$ and $\Phi$ are finite, $\fhi^0(t,\cdot)\in V_0$,
and $\vu(t,\cdot)\in \Vn$, then $\vu$ is divergence-free
and, by \eqref{zeromat}, $\frac{D\fhi}{Dt}$ has zero mean value.

%
%
%
%
%
We also note that, whenever $\frac{D\fhi}{Dt} \in H_0$, then,
by elliptic regularity, it follows
\[
  \mu^0\in V_0\cap H^2(\Omega)
   \quad\text{and}\ \  \nabla_x \mu^0 \cdot {\bf n}_{|_{\Gamma}}=0.
\]
Moreover, we can equivalently rewrite the pseudopotential $\Phi$ as follows:
\begin{equation}\label{phi0}
  \Phi(\delta E, E)=\int_\Omega\widetilde\phi(\delta E, E)\, \dix,
  \quad \hbox{where}\ \ \widetilde\phi(\delta E, E)
  = \phi(\delta E, E)+\frac12|\nabla_x\mu^0|^2.
\end{equation}
Indeed, integrating by parts in space and using the definition of $J$
\eqref{duale}, it turns out that:
\begin{equation}\label{phi0b}
  \io|\nabla_x\mu^0|^2\,\dix=-\io\Delta\mu^0\mu^0\,\dix
   =\io J(\mu^0)\mu^0\,\dix=\left\langle\frac{D\fhi}{Dt} , J^{-1}\left(\frac{D\fhi}{Dt} \right)\right\rangle.
\end{equation}
An alternative strategy to derive the (isothermal) Cahn-Hilliard equation
starting from a balance of the so-called ``microscopic motions''
can be found in \cite{Gurtin} and could be extended to the present
case of binary fluids. However, while the analysis of the isothermal Cahn-Hilliard equation
has been deeply investigated starting from the pioneering paper \cite{cahnhill} 
and up to the most recent contributions and Navier-Stokes systems 
(cf., e.g., \cite{A09, GG1, GPV96}), at least up to our knowledge
a rigorous derivation of a thermodynamically consistent model
including also the internal energy balance equation 
is still lacking (cf., e.g. \cite{rocca} for a thermodynamically
consistent model coupling the Cahn-Hilliard system with a singular heat equation).
Our method is intended to fill this gap in the more intricated
case where the effects of the macroscopic velocity
$\bu$ are also taken into account.


\subsection{Constitutive relations}

We start by introducing the energy density
$B$ and the energy flux vector ${\bf H}$, both assumed to be the sum of their
non-dissipative and dissipative components, namely, $B=B^{nd}+B^d$,
${\bf H}={\bf H}^{nd}+{\bf H}^d$,
where
\begin{gather}\label{Bnd}
  B^{nd}
   =\frac{\partial \psi }{\partial \varphi}
   = \frac1\varepsilon F'(\varphi) - \theta,\\
 \label{Bd}
  B^{d}=\delta_{H_0,\frac{D\varphi}{Dt}} \Phi
    =  J^{-1}\left(\frac{D\varphi}{Dt}\right),\\
 \label{Hnd}
  {\bf H}^{nd}=\frac{\partial \psi }{\partial \Grad \varphi} = \varepsilon \Grad \varphi.
\end{gather}
Moreover, we set  ${\bf H}^d \equiv 0$. Here and in the sequel for simplicity 
we will use the symbol $\partial$ not only for partial derivatives, but also 
for variational derivatives of possibly non-convex funtionals. 
Relation \eqref{Bd} defines $B^d$
as the first variation (more precisely, subdifferential)
of $\Phi$ with respect to $\frac{D\varphi}{Dt}$ in the space $H_0$. To see that it
coincides, indeed, with $J^{-1}\left(\frac{D\varphi}{Dt}\right)$, we observe
that, for any $v\in H_0$, there holds
\begin{align}\nonumber
  & \left( J^{-1}\left(\frac{D\varphi}{Dt}\right), v - \frac{D\varphi}{Dt} \right)
    = \left\langle v - \frac{D\varphi}{Dt}, J^{-1}\left(\frac{D\varphi}{Dt}\right) \right\rangle
    = \bigg(\!\!\bigg( v - \frac{D\varphi}{Dt}, \frac{D\varphi}{Dt} \bigg)\!\!\bigg)_{V_0'}\\
 \label{vprim}
  & \mbox{}~~~~~
  \le \frac12 \| v \|_{V_0'}^2 - \frac12 \Big\| \frac{D\varphi}{Dt} \Big\|_{V_0'}^2
   = \frac12 \langle v, J^{-1} v \rangle
    - \frac12 \left\langle \frac{D\varphi}{Dt}, J^{-1}\left(\frac{D\varphi}{Dt}\right) \right\rangle.
\end{align}
%

The dissipative components of the
heat and entropy fluxes (denoted respectively
by ${\bf q}^d$ and ${\bf Q}^d$) are
\begin{equation}\label{Q}
  {\bf q}^d=\theta {\bf Q}^d =
   -\theta \frac{\partial \phi }{\partial \Grad\theta }%
   =-\kappa(\theta )\Grad\theta,
\end{equation}
whereas the nondissipative components
${\bf q}^{nd}$ and ${\bf Q}^{nd}$ will be determined later
on in such a way that the second law of Thermodynamics
is satisfied (cf.~\eqref{clau} below).
Of course, we assume that
${\bf q}={\bf q}^{d}+{\bf q}^{nd}$,
${\bf Q}={\bf Q}^{d}+{\bf Q}^{nd}$,
and ${\bf q}^{nd}=\theta {\bf Q}^{nd}$.
In what follows we also ask that 
$\kappa(\theta)\sim 1+\theta^\beta$ with $\beta\geq 2$ (cf.~\eqref{hpbetadelta} below).
This choice is mainly motivated by mathematical reasons 
(actually, it guarantees some additional 
integrability of $\teta$). A physical justification for it is provided, e.g., in \cite{zr}.

Also the stress tensor $\sigma$ is decomposed into the dissipative
component
\begin{equation}\label{sigmad}
  \sigma ^{d}=\frac{\partial \phi }{\partial D{\bf u}}
   =\nu (\theta )D{\bf u} - p\mathbb{I}=:\tn{S}-p\tn{I},
\end{equation}
\[
  - p \in \partial I_0 (\dive \bu), \ \tn{S} = \nu(\theta) D{\bf u},
\]
and the non dissipative part $\sigma^{nd}$ to be
determined below.

The entropy of the system is given by
\begin{equation}\label{s}
  s = -\frac{\partial \psi }{\partial \theta }=-f'(\teta)+\varphi,
\end{equation}
and, finally, the internal energy $e$ reads
\begin{equation}\label{e}
  e =\psi +\theta s= \frac1\varepsilon F(\varphi)
    + \frac{\varepsilon}{2} {|\Grad\varphi|^2}+Q(\teta),
\end{equation}
where $Q(\teta)=f(\teta)-\teta f'(\teta)$
represents, physically speaking, the 
antiderivative of the specific heat $c_V$.


\subsection{Field equations}
\label{sec:pdes}

In accordance with \emph{Newton's second law}, the balance of momentum reads
\begin{equation}\label{mombal}
  \partial_t \ub  + \dive (\ub \otimes \ub )  = \dive \sigma + \vc{g},
\end{equation}
where $\vc{g}$ is a given external force.

The balance of internal energy takes the form
\begin{equation}\label{eqe}
 \frac{D e}{D t}+\dive {\bf q}
  =\sigma :D(\ub)
 + B \frac{D\varphi}{Dt} + {\bf H} \cdot \Grad\frac{D\varphi}{Dt}+{N},
\end{equation}
where we recall that the internal energy flux
is decomposed as $\vc{q} = \vc{q}^d + \vc{q}^{nd}$.
Moreover, we notice that, on the right hand side of \eqref{eqe},
there appears a new (with respect to the standard theory of \cite{Frem})
term ${ N}$. This contribution is aimed at balancing the
nonlocal dependence of the last term in the pseudopotential of dissipation $\Phi$
(cf.~\eqref{Phi}) with respect to the dissipative variable $\frac{D\fhi}{D t}$, 
which is  one of the peculiarities of Cahn-Hilliard systems (cf.~also \cite{KRS07}
for similar techniques applied to different nonlocal contributions).
The expression of $N$ (in terms of the dissipative variables $\delta E$)
will be obtained below in such a way to comply with
the second law of Thermodynamics. Since the main role of ${N}$ is to model
the nonlocal interactions between points inside $\Omega$, it will result
from our computations that $\io {N}(x)\,\dix=0$, in agreement with natural expectations.

Finally, the equation ruling the evolution of the order
parameter $\varphi$ can be derived from the
principle of virtual powers. Indeed, following the general theory
developed in~\cite[Chap. 2]{Frem}, we have
\begin{equation} \label{pom}
  \dive {\bf H} - B = 0.
\end{equation}
However, in the present setting the above relation is not
completely rigorous, since it does not properly incorporate the
boundary conditions and the mass conservation constraint. So, to be more precise, we
%
%
first have to rewrite the expression of $\Psi$ in the form
\begin{equation}\label{psi0}
  \Psi(E) = \io \left(\frac{\varepsilon}{2}|\Grad\varphi^0|^{2}
  + \frac1\varepsilon F(\fhi^0+m_0)
  + f(\theta)-\teta(\varphi^0+m_0)\right) \, \dix,
\end{equation}
where we used the decomposition $\fhi = \fhi^0 + m_0$ already introduced in Section \ref{free-energy-sec}.

In order to impose this constraint mathematically, we restate \eqref{pom} as
a generalized gradient flow problem {\em in the space $H_0$}\/
as follows:
\begin{equation} \label{pom2}
  B^d + \delta_{H_0,\varphi^0} \Psi
   = \delta_{H_0,\frac{D\varphi}{Dt}} \Phi
    + \delta_{H_0,\varphi^0} \Psi = 0.
\end{equation}
Let us note that asking for $\fhi^0$ to lie in the domain of the
differential $\delta_{H_0,\varphi^0} \Psi$
means that there exists a (unique) function  $z \in H_0$
such that $\delta_{H_0,\varphi^0} \Psi(\fhi^0)$
can be represented by $z$ in the scalar product of $H_0$
(i.e., of $H$).
In this way, \eqref{pom2} incorporates both the homogeneous Neumann boundary
conditions for $\fhi$ and the mass conservation property.

Moreover, it is immediately seen that such a function
$z$ must have the expression
\begin{equation} \label{eqd0}
  z = - \varepsilon \Delta \varphi^0
   + \frac1\varepsilon \left( F'(\varphi^0+m_0) - \overline{F'(\varphi^0+m_0)} \right)
   - \theta + \overline{\theta}.
\end{equation}
Combining \eqref{pom2} and \eqref{eqd0} with \eqref{Bd}, we then get
\begin{equation} \label{eqd1}
  J^{-1} \left( \varphi_t + \ub \cdot \Grad \varphi \right)
   = \varepsilon \Delta \varphi^0
   - \frac1\varepsilon \left( F'(\varphi^0+m_0) - \overline{F'(\varphi^0+m_0)} \right)
   + \theta - \overline{\theta}.
\end{equation}
Applying the distributional Laplace operator to both hand sides
and noting that $- \Delta J^{-1} v = v$ for any $v \in H_0$
(cf.~\eqref{duale}), we then arrive at the system
\begin{align} \label{eqd2}
  & \varphi_t + \ub \cdot \Grad \varphi = \Delta \mu,\\
 \label{eqd3}
  & \mu = - \varepsilon \Delta \varphi + \frac1\varepsilon F'(\varphi) - \teta,
\end{align}
where the auxiliary variable $\mu$ is introduced mainly for
mathematical convenience and takes the name of {\em chemical potential}.
In fact, we may note that  $\mu^0 =  -J^{-1}(\frac{D\varphi}{Dt})=\mu - \overline{\mu}$ 
(cf.~\eqref{eqd1}, \eqref{muzero}).

The non-dissipative components of the stress $\sigma^{nd}$ and of the flux
$\vc{q}^{nd}$, as well as the ``nonlocality compensation'' term $N$,
are determined by means of \eqref{eqe} and the constitutive relations derived above,
in order for the {\it second law of Thermodynamics}\/ to be satisfied.
Indeed, computing $\frac{D e}{D t}$ from \eqref{e}
by means of the standard Helmholtz relations 
\begin{equation}\label{Helmholtz}
e=\psi+\theta s, \quad s=-\frac{\partial\psi}{\partial\theta},
\end{equation}
we get
\begin{equation}\label{ide}
  \frac{D e}{D t}
   = \frac{D\psi}{ Dt}
    + \theta\frac{D s}{ D t}
    + \frac{D \theta}{D t} s
   = \frac{\partial \psi}{\partial\varphi} \frac{D\varphi}{D t}
    + \frac{\partial\psi}{\partial \Grad\varphi} \cdot \frac{D(\Grad\varphi)}{D t}
    + \theta\frac{D s}{D t},
\end{equation}
whereas
\begin{equation}\label{psigradb}
  \frac{\partial\psi}{\partial\Grad\varphi} \cdot \frac{D (\Grad\varphi)}{D t}
   = {\bf H}^{nd} \cdot \left( \Grad\frac{D \varphi}{D t}
   - \Grad\ub \cdot \Grad\varphi\right).
\end{equation}
Moreover, by \eqref{eqd0}--\eqref{eqd3}, we have (pointwise)
\begin{equation}\label{theta2}
  \frac{D\varphi}{Dt} J^{-1}\left(\frac{D\varphi}{Dt}\right)
   = - \Delta \mu ( \mu - \overline{\mu} )
   = - \frac12 \Delta (\mu - \overline{\mu})^2 + | \nabla_x \mu |^2.
\end{equation}
To deduce the expressions for the non-dissipative
components of the stress $\sigma^{nd}$ and of the flux
$\vc{q}^{nd}$ as well as that of $N$, we impose validity
of the Clausius-Duhem inequality in the form
\begin{equation}\label{clau}
  \teta\Big( \frac{D s}{Dt} + \dive {\bf Q} \Big) \geq 0,
\end{equation}
where ${\bf Q}$ denotes the entropy flux and it is linked to the flux ${\bf q}$ by the relation
$\theta{\bf Q}={\bf q}$.
Recalling \eqref{Bnd}--\eqref{Q} and
developing the left hand side, we get
%
%
\begin{align}\nonumber
  \teta\left(\frac{D s}{Dt} + \dive {\bf Q}\right)
   \overset{\eqref{e}}{=} &\frac{D e}{Dt}+\dive {\bf q}-\frac{D\psi}{Dt}
    -\frac{D\teta}{Dt}s-{\bf Q}\cdot\nabla_x\teta\\
 \nonumber
  \overset{\eqref{Helmholtz}}{=} &\frac{D e}{Dt}
    + \dive {\bf q}- \frac{\partial\psi}{\partial\varphi}\frac{D\varphi}{Dt}
    - \frac{\partial\psi}{\partial\nabla_x\varphi} \cdot \frac{D\nabla_x\varphi}{Dt}
    -{\bf Q}\cdot\nabla_x\teta\\
 \nonumber
 \overset{\eqref{eqe}, \, \eqref{Bnd}}{=} &\sigma:D(\bu)+B\frac{D\varphi}{Dt}
   + {\bf H}\cdot\frac{D\nabla_x\varphi}{Dt} + N
   - \frac{\partial\psi}{\partial\nabla_x\varphi}  \cdot \frac{D\nabla_x\varphi}{Dt}\\
 \nonumber
  & \mbox{}~~~~~
   - B^{nd} \frac{D\varphi}{Dt}
   + \frac{\kappa(\teta)}{\teta}|\nabla_x\teta|^2
   - {\bf Q}^{nd}\cdot\nabla_x\teta\\
 \nonumber
  \overset{\eqref{Bd}, \eqref{psigradb},\,\eqref{theta2},\, \eqref{sigmad}}{=} &
    \big(\sigma^{nd}+\varepsilon(\nabla_x\varphi\otimes\nabla_x\varphi)\big):D(\bu)
   + \nu(\theta)|D(\bu)|^2 +|\nabla_x\mu|^2
   - \frac12\Delta(\mu-\overline{\mu})^2+N\\
\nonumber
  & \mbox{}~~~~~
   + \frac{\kappa(\teta)}{\teta}|\nabla_x\teta|^2
   - {\bf Q}^{nd}\cdot\nabla_x\teta.
\end{align}
Then, in order to obtain the non-negativity of the right hand
side  (cf.~\eqref{clau}), we can assume, e.g., the following constitutive relations 
\begin{equation}\label{sigmand}
  \sigma^{nd} = -\varepsilon \Grad \varphi\otimes\Grad\varphi, \quad
   \vc{q}^{nd} = 0, \quad N=\frac12\Delta(\mu-\overline{\mu})^2.
\end{equation}
With these choices, we get $\io N(x)\, \dix=0$, as expected. Moreover,
a straighforward computation shows that
the internal energy balance \eqref{eqe}
can be rewritten as
\begin{equation}\label{ebal0}
  (Q(\theta))_t + {\bf u} \cdot \nabla_x Q(\theta) +\theta \frac{D \varphi}{D t}
  - \dive (\kappa(\theta) \nabla_x \theta) = \nu(\theta) |D {\bf u}|^2 + |\nabla_x \mu|^2.
\end{equation}
Notice that the dissipation terms on the right hand side are in perfect
agreement with the expression \eqref{phi0} of the pseudopotential
of dissipation $\Phi$. Indeed, as already mentioned,  one has $\mu^0=\mu-\overline{\mu}$
due to \eqref{eqd2}--\eqref{eqd3}.


\subsection{Strong formulation}
\label{subs:strong}

On account of the derivation sketched above, we can now write the PDE
system representing the {\em strong formulation}\/ of our model.

Firstly, collecting \eqref{sigmad}, \eqref{mombal} and
\eqref{sigmand}, we obtain that the evolution of the velocity
$\vu$ is ruled by the Navier-Stokes system, given by

\medskip

\greybox{
\centerline{\textsc{Incompressibility:}}

\vspace{-2mm}

\begin{equation} \label{incom}
  \dive \bu = 0;
\end{equation}

\smallskip

\centerline{\textsc{Conservation of momentum:}}

\vspace{-2mm}

\begin{equation} \label{momentumbal}
  \bu_t + \bu \cdot \Grad \bu + \Grad p = \dive \tn{S}
  - \dive \left( \varepsilon \Grad \varphi \otimes\Grad\varphi\right)  + \vc{g},
\end{equation}
}
\noindent%
where $p$ is the pressure, and
\begin{equation}\label{stress}
  \mathbb{S} =  \frac{\nu(\theta)}{2} \left( \Grad \bu + \Grad^t \bu \right).
\end{equation}

Regarding the evolution of $\varphi$, \eqref{eqd2}--\eqref{eqd3}
give rise to the following

\medskip

\greybox{
\centerline{\textsc{Cahn-Hilliard equation:}}

\vspace{-2mm}

\begin{equation}
\label{CH}
{\textnormal{
$\left \{
\begin{array}{lll}
\!\!\!\!\!\! & \varphi_t + {\bf u} \cdot \nabla_x \varphi = \Delta \mu, \\[2mm]
\!\!\!\!\!\! & \mu = - \varepsilon\Delta \varphi
   + \frac{1}{\varepsilon} F'(\varphi)  - \theta.
\end{array}
\right.$}}
\end{equation}
}
\noindent%
Finally, we rewrite \eqref{ebal0} as the

\medskip

\greybox{
\centerline{\textsc{Internal energy balance:}}

\vspace{-2mm}

\begin{equation}\label{theta}
  (Q(\theta))_t + {\bf u} \cdot \nabla_x Q(\theta) +\theta \frac{D \varphi}{D t}
    - \dive (\kappa(\theta) \nabla_x \theta) = \nu(\theta) |D {\bf u}|^2
    + |\nabla_x \mu|^2.
\end{equation}
}
\noindent%
It is worth observing that, neglecting the temperature,
equations (\ref{incom}--\ref{theta})
reduce to the model derived in \cite{GPV96}
by means of a different method.


\subsection{Balances for total energy and entropy}
\label{subs:enen}

A key point in the statement of the weak formulation
of our model consists in replacing the
``heat'' equation \eqref{theta} with the balances of
total energy and of entropy. These relations are, indeed,
mathematically more tractable, but they keep all the main 
features of the problem (in particular,
the first and second laws of Thermodynamics are still respected). 

We give here a formal derivation of these relations, beginning
with the total energy balance.
To obtain it, we start multiplying \eqref{momentumbal} by $\vu$,
which gives
\begin{align} \nonumber
  & \frac12 \ddt |\bu|^2 + \frac12 \bu \cdot \Grad | \bu |^2 + \dive ( p \vu )
  = \dive ( \tn{S} \vu ) - \tn{S} : \Grad \vu
  - \varepsilon \dive \big( ( \Grad \varphi \otimes \Grad\varphi) \vu \big) \\
 \label{en11}
 & \mbox{}~~~~~
  + \varepsilon \left( \Grad \varphi \otimes \Grad\varphi \right) : D \vu
  + \vc{g} \cdot \vu.
\end{align}
Next, we multiply the first \eqref{CH} by $\mu$, the second by $\fhi_t$, and
take the difference. After standard manipulations, we obtain
\begin{equation} \label{en12}
  \ddt \left( \varepsilon \frac{|\Grad \fhi|^2}2 + \frac1{\varepsilon} F(\fhi) \right)
   - \theta \fhi_t - \varepsilon \dive (\fhi_t \Grad \fhi)
   - \dive ( \mu \Grad \mu ) + | \Grad \mu |^2 + \vu \cdot \Grad \fhi \mu
  = 0.
\end{equation}
We now substitute into the last term the expression of $\mu$
given by the second \eqref{CH}:
\begin{align} \nonumber
  \vu \cdot \Grad \fhi \mu
   & = \frac{1}{\varepsilon} \vu \cdot \Grad F (\fhi)
    - \theta \vu \cdot \Grad \fhi
    - \varepsilon \vu \cdot \Grad \fhi \Delta \fhi \\
 \nonumber
   & = \frac{1}{\varepsilon} \vu \cdot \Grad F (\fhi)
    - \theta \vu \cdot \Grad \fhi
    - \varepsilon \dive \big( ( \Grad \varphi \otimes \Grad\varphi ) \vu \big)\\
 \label{en13}
   & \mbox{}~~~~~~~~~~~~~~~~~~~~~
    + \varepsilon \left( \Grad \varphi \otimes \Grad\varphi \right) : D \vu
    + \varepsilon \vu \cdot \Grad \frac{ | \Grad \fhi |^2 }2.
\end{align}
By \eqref{en13}, \eqref{en12} is transformed into
\begin{align} \nonumber
  & \ddt \left( \varepsilon \frac{|\Grad \fhi|^2}2 + \frac1{\varepsilon} F(\fhi) \right)
   - \theta \fhi_t - \varepsilon \dive (\fhi_t \Grad \fhi)
   - \dive (\mu \Grad \mu) + | \Grad \mu |^2
   + \frac{1}{\varepsilon} \vu \cdot \Grad F (\fhi) \\
  \nonumber
   & \mbox{}~~~~~~~~~~
   - \theta \vu \cdot \Grad \fhi
   - \varepsilon \dive \big( ( \Grad \varphi \otimes \Grad\varphi ) \vu \big)
   + \varepsilon \left( \Grad \varphi \otimes \Grad\varphi \right) : D \vu
   + \varepsilon \vu \cdot \Grad \frac{ | \Grad \fhi |^2 }2 \\
 \label{en14}
   & \mbox{}~~~~~ = 0.
\end{align}
Then, we can take the sum of \eqref{theta}, \eqref{en11} and \eqref{en14}
to obtain the

\medskip

\greybox{
\centerline{\textsc{Total energy balance:}}

\vspace{-2mm}

\begin{equation} \label{energy0}
  \partial_t \left( \frac{1}{2}|\bu|^2 + e \right)
   + \bu \cdot \Grad \left( \frac{1}{2}|\bu|^2 + e \right)
   + \dive \Big( p \bu + \vc{q} - \tn{S} \bu\Big)
   - \dive \big ( \varepsilon \fhi_t \Grad \fhi + \mu \Grad \mu )
  = \vc{g} \cdot \bu,
\end{equation}
}
\noindent%
with the internal energy
\begin{equation}\label{int_pot_en}
  e = \frac1\varepsilon F(\varphi)+\frac{\varepsilon}{2}|\Grad \varphi |^2  + Q(\theta)
\end{equation}
and the heat flux
\begin{equation}\label{hflux}
   {\bf q}
  = -\kappa( \theta )\Grad\theta.
\end{equation}
Next, multiplying \eqref{theta} by $1/\theta$, we obtain the

\medskip

\greybox{
\centerline{\textsc{Entropy equation:}}

\vspace{-2mm}

\begin{equation}\label{entropy}\displaystyle
  (\Lambda(\theta) + \varphi)_t
   + {\bf u} \cdot \nabla_x (\Lambda(\theta)) + {\bf u} \cdot \nabla_x \varphi
   - \dive \left (\frac{\kappa(\theta) \nabla \theta}{\theta} \right )
  = \frac{\nu(\theta)}{\theta} |D {\bf u}|^2
   + \frac{1}{\theta} |\nabla_x \mu|^2
   + \frac{\kappa(\theta)}{\theta^2} |\nabla_x \theta|^2,
\end{equation}
}
\noindent%
where
\begin{equation}
\label{Lambda}
\Lambda(\theta) = \int_1^{\theta} \frac{c_V(s)}{s} \, \dis.
\end{equation}
We anticipate that the function $\Lambda$ is well defined in
view of the assumptions on $c_V$ stated in \eqref{kappacv} below.
Moreover, it is worth noting that \eqref{entropy}
is an equality at this level, but it will turn to an inequality
(cf.~\eqref{wf6} below) in the framework of
the rigorous definition of \emph{weak solution}\/ that
will be introduced later on (cf.~Def.~\ref{defweak}
and Theorem~\ref{thm:main}).
Of course, this phenomenon is due to the
quadratic terms on the right hand side, which do not behave
well with respect to weak limits.
However, it is worth noticing that, in case we could prove
that there exist a {\em smooth} solution to the weak
formulation of the model, then for that solution
it would be possible to recover the
entropy {\em equality}~\eqref{entropy}. Moreover, \eqref{entropy}
is equivalent to \eqref{eqe} in that setting. In this
sense, the weak formulation analyzed in the next
section turns out to be compatible with the strong formulation
(\ref{model1}--\ref{model5}) given at the beginning,
at least when sufficient regularity holds.


\subsection{Initial and boundary conditions}
\label{IBC}

In order to get a well-posed problem, we have to specify
suitable initial and boundary conditions.
Compatibly with the physical derivation, we will
essentially assume that
the system is insulated from the exterior.
This leads to taking the following no mass flux
(through the boundary) condition:
\begin{equation}
\label{BC2}
\nabla_x \mu \cdot {\bf n}_{|_{\Gamma}} = 0,
\end{equation}
where we recall that ${\bf n}$   is the external normal.
Next, we assume that
\begin{equation}
\label{BC3}
\nabla_x \varphi \cdot {\bf n}_{|_{\Gamma}} = 0.
\end{equation}
This position prescribes a ``contact angle'' of $\pi/2$ between the diffuse
interface and the boundary of the domain.
Moreover, we take no-flux boundary conditions for the temperature:
\begin{equation}
\label{BC4}
{\bf q} \cdot {\bf n}_{|_{\Gamma}} = 0.
\end{equation}
Finally, we assume {\it complete slip} boundary condition in the momentum
equation \eqref{momentumbal}:
\begin{equation}
\label{BC1}
  {\bf u} \cdot {\bf n}_{|_{\Gamma}} = 0,
   \qquad [ \mathbb{S} {\bf n}] \times {\bf n}_{|_{\Gamma}} = 0.
\end{equation}
The first condition states that the normal component of the boundary
velocity is zero (so, the fluid cannot exit from $\Omega$,
but it can move tangentially to the boundary).
The second position prescribes that there is no external contribution
to the viscous stress.
This, in a sense, excludes friction effects with the boundary. The
above choice has also mathematical implications. Indeed,
following the lines of \cite{BFM09}, we will use it in order
estimate the pressure term appearing in the total energy
balance \eqref{energy0}.
It is worth noting that an estimate of the pressure can be reached also
in the case when $\Omega$ is the unit torus and periodic boundary conditions are
taken for all unknowns. In particular, our results could be extended to that
situation with trivial modifications. For other types
of boundary conditions, integrability of
the pressure is, instead, an open issue.

We remark once more that, thanks to \eqref{BC2}
and the first \eqref{BC1},
integrating \eqref{CH} in space and time, we
get mass conservation:
\begin{equation}\label{massCons}
  \io\varphi(t)=\io\varphi(0)\quad\forall\, t\in (0,T),
\end{equation}
This feature is characteristic of Cahn-Hilliard-type models and is in
agreement with the underlying physics.

Finally, the system is complemented by the initial conditions
\begin{equation}\label{IC}
  {\bf u}(0, \cdot) = {\bf u}_0, \qquad 
   \varphi(0, \cdot) = \varphi_0, \qquad 
   \theta(0, \cdot) = \theta_0.
\end{equation}
%


\section{Main result}
\label{mainres}


\subsection{Assumptions on coefficients and data}

\label{Sect-Hp}

Before formulating the main result of the paper, we list here the
hypotheses imposed on the constitutive functions.
First of all, just for the sake of simplicity, we take
$\varepsilon=1$ and ${\bf g}={\bf 0}$. The case
of a nonconstant forcing term could be treated, indeed,
with trivial modifications.
Next, we assume that $F(\varphi)$ is the classical double-well potential, namely
\[
F(\varphi) = \frac{1}{4} (\varphi^2 - 1)^2,
\]
so that
\begin{equation}\label{effep}
F'(\varphi) = \varphi^3 - \varphi.
\end{equation}
 More general expressions of $F'$ having cubic growth at $\infty$ may be admissible as well, but we prefer to keep
right from the beginning the expression \eqref{effep} in order not to overburden the presentation.
We assume that the thermal conductivity, the specific heat {and the viscosity of the mixture}
depend on $\theta$ in the following way:
\begin{equation}
\label{kappacv}
  \kappa(\theta) = 1 + \theta^{\beta},
   \qquad c_V(\theta) = \theta^{\delta},
   \qquad 0 < \underline{\nu} \le \nu(\theta) \le \overline{\nu},
\end{equation}
for all $\theta\geq 0$, some $0<\underline{\nu} <\overline{\nu}$,
and some $\beta > 0, \delta > 0$
complying with the following restrictions:
\begin{equation}
\label{hpbetadelta}
{\beta \geq 2}, \qquad \frac{1}{2} < \delta < 1.
\end{equation}
In view of \eqref{kappacv}, we can compute
\begin{equation}\label{defQ}
  Q(\theta) = \int_0^\theta c_V(s) \, \dis
   = \frac1{\delta + 1} \theta^{\delta+1},
\end{equation}
as well as (cf.~\eqref{Lambda})
\begin{equation}\label{Lam}
  \Lambda(\theta)
   = \frac1{\delta} (\theta^{\delta} - 1 ).
\end{equation}
The ansatz $\delta < 1$ comes from the assumption of physical
consistency (actually, it implies that
the thermal component $\Lambda$ of the entropy is concave,
as prescribed by Thermodynamics);
the other limitations mainly have
a mathematical motivation and are needed in order to obtain
the necessary regularity to pass to the limit. Notice however that
a power-like behavior for the heat-conductivity is typical of several
types of fluids (cf., e.g., \cite{zr}).

For brevity we also set
\begin{equation}
\label{cost_A}
{p_{\beta, \delta}} := \beta + \frac{2}{3} (\delta + 1).
\end{equation}
This exponent will be needed in the a-priori estimates
derived below (cf.~for instance \eqref{XY} and \eqref{conv_theta}).

We conclude by specifying our hypotheses on the initial data:
\begin{equation}\label{hyp:0}
  \vu_0 \in L^2_{\dive}(\Omega;\RR^3), \qquad
   \fhi_0 \in H^1(\Omega), \qquad
   \theta_0\in L^{\delta+1}(\Omega),~~~\theta_0>0~\text{almost everywhere.}
\end{equation}
Here and below, $L^2_{\dive}$ indicates the space of divergence-free $L^2$
functions.


\subsection{Weak formulation}

\label{sect-weak}

First of all, we rewrite the momentum equation \eqref{mombal},
with ${\bf g}={\bf 0}$, in the more explicit form:
\begin{equation}
\label{momentum}
  {\bf u}_t + {\bf u} \cdot \nabla_x {\bf u} + \nabla_x p
  = \dive \mathbb{S}
  - \dive ( \Grad \fhi \otimes \Grad \fhi).
\end{equation}
%
%
This permits us to introduce the notion of {\it weak solution} to our model problem:
\begin{defn}\label{defweak}
A {\rm weak solution} to the non-isothermal diffuse interface model
for two-phase flows of fluids is a quadruplet $({\bf u}, \varphi, \mu, \theta)$
satisfying the {\bf incompressibility}\/ condition
$\dive \vu = 0$ a.e.~in~$(0,T)\times\Omega$,
the {\bf weak momentum balance}
\begin{align} \nonumber
  & \int_0^T \int_{\Omega} ({\bf u} \cdot \partial_t \xibf
  + ({\bf u} \otimes {\bf u}) : \nabla_x \phibf + p \, \dive \xibf) \\
 \label{wf1}
  & = \int_0^T \int_{\Omega} (\mathbb{S} : \nabla_x \xibf)
  - \int_0^T \int_{\Omega} ( \Grad \fhi \otimes \Grad \fhi ) : \nabla_x \xibf
  - \int_{\Omega} \bu_0 \cdot \xibf(0, \cdot),
\end{align}
for all $\xibf \in \mathcal{C}^{\infty}_0([0,T) \times \overline{\Omega}; \mathbb{R}^3)$
such that ${\xibf \cdot {\bf n}}_{|_{\Gamma}} = 0$, the
{\bf Cahn-Hilliard system}
\begin{align}\label{wf3}
  & \langle \varphi_t, \xi \rangle + \io ( {\bf u} \cdot \nabla_x \varphi ) \xi
   = \io \Grad \mu \cdot \Grad \xi \qquad \textnormal{for all $\xi \in V$, and a.e.~in $(0,T)$},\\
 \label{wf4}
  & \mu = - \Delta \varphi + F'(\varphi) - \theta \qquad \textnormal{a.e. in $(0,T)\times \Omega$},
\end{align}
with the boundary condition \eqref{BC3} and the initial condition
${\fhi}(0, \cdot) = \fhi_0$, and the {\bf weak total energy balance}
\begin{align}
& \!\!\!\!\!\!\! \int_0^T \int_{\Omega} \left ( \frac{1}{2} |{\bf u}|^2 + e \right ) \partial_t \xi
  + \int_0^T \int_{\Omega} \left ( \frac{1}{2} |{\bf u}|^2 \, {\bf u}
         + e \, {\bf u} \right ) \cdot \nabla_x \xi
  + \int_0^T\int_\Omega \hat{\kappa}(\theta) \Delta\xi +\int_0^T\io p{\bf u}\cdot\nabla\xi \nonumber\\
 & - \int_0^T \int_{\Omega} ( \mathbb{S} \, {\bf u})\cdot \nabla_x \xi
  + \int_0^T \int_{\Omega} \frac{\mu^2}{2} \Delta \xi
  + \int_0^T \int_{\Omega} ({\bf u} \cdot \nabla_x \varphi) \, (\nabla_x \varphi \cdot \nabla_x \xi)
   + \int_0^T \int_{\Omega} (\nabla_x \mu \otimes \nabla_x \xi) : \nabla_x \nabla_x \varphi \nonumber\\
 & 
  + \int_0^T \int_{\Omega} (\nabla_x \mu \otimes \nabla_x \varphi) : \nabla_x \nabla_x \xi 
   - \int_{\Omega} \left (\frac{1}{2} |{\bf u}_0|^2 + e_0 \right ) \, \xi(0, \cdot) 
    = 0,\quad\hbox{for all $\xi\in \mathcal{C}^{\infty}_0([0,T) \times \Omega)$},
 \label{wf5}
\end{align}
where $e$ is given by \eqref{int_pot_en}, $\hat{\kappa}$ is defined as
\begin{equation}\label{kciapo}
  \hat{\kappa}(\theta) = \int_0^\theta \kappa(s) \, \dis
   = \theta + \frac1{\beta + 1} \theta^{\beta+1},
\end{equation}
and finally, we have set $e_0 = F(\varphi_0) + \frac{|\nabla_x \varphi_0|^2}{2} + Q(\theta_0)$.
%
\end{defn}
\noindent%
It is worth noting that \eqref{wf1}
incorporates both the incompressibility constraint and
the initial condition ${\bf u}(0, \cdot) = {\bf u}_0$; moreover,
it accounts for the complete-slip conditions \eqref{BC1}.
The first equation \eqref{wf3} of the Cahn-Hilliard system is in
weak form and also accounts for the
no-flux condition~\eqref{BC2},
while we will be able to prove sufficient regularity on $\fhi$ in
order for \eqref{wf4} to hold pointwise (with the no-flux condition
\eqref{BC3} in the sense of traces).
To get \eqref{wf5}, we  tested \eqref{energy0}  
by $\xi$, integrated by parts in time and used
the Cahn-Hilliard system \eqref{CH}. More precisely, we wrote
\begin{align}\label{towf11}
  & - \io \fhi_t \Grad \fhi \cdot \Grad \xi
   = \io ( \vu \cdot \Grad \fhi ) \, ( \Grad \fhi \cdot \Grad \xi )
   - \io \Delta \mu \, ( \Grad \fhi \cdot \Grad \xi ) \\
 \nonumber
  & \mbox{}~~~~~
   = \io ( \vu \cdot \Grad \fhi ) \, ( \Grad \fhi \cdot \Grad \xi )
   + \io ( \Grad \mu \otimes \Grad \xi ) : \Grad \Grad \fhi
   + \io ( \Grad \mu \otimes \Grad \fhi ) : \Grad \Grad \xi.
\end{align}
%
%


\subsection{Main existence theorem}
\label{subsec:main}

Our main result reads as follows:
\begin{thm}\label{thm:main}
Under the assumptions stated in\/ {\rm Subsection~\ref{Sect-Hp}}, the non-isothermal
diffuse interface model for two-phase flows of fluids admits at least a weak solution,
in the sense of\/ {\rm Definition~\ref{defweak}}, in the following regularity class:
\begin{eqnarray}
&& {\bf u} \in L^{\infty}(0,T; L^2(\Omega; \mathbb{R}^3)) \cap L^2(0,T;\Vn)  \label{reg1}\\
&& \varphi \in H^1(0,T; (H^1(\Omega))')\cap L^{\infty}(0,T; H^1(\Omega)) \cap L^2(0,T; H^3(\Omega)) \label{reg2}\\
&& \mu \in L^2(0,T; H^1(\Omega)) \cap L^{\frac{14}{5}}((0,T) \times \Omega) \label{reg3}\\
&& \theta \in L^{\infty}(0,T; L^{\delta+1}(\Omega)) \cap L^{\beta}(0,T; L^{3 \beta}(\Omega))
 \cap L^2(0,T; H^1(\Omega))\label{reg4}\\
&& \theta > 0~~\text{a.e.~in $(0,T)\times \Omega$},\quad
  \log\theta\in L^2(0,T; H^1(\Omega)), \label{reg4b}
\end{eqnarray}
$\delta$ and $\beta$ being specified in \eqref{hpbetadelta}.
Moreover, this solution complies with the following weak form of the
{\bf entropy production inequality}:
\begin{align} \nonumber
 & \int_0^T\int_{\Omega} (\Lambda(\theta) + \varphi) \, \partial_t \xi
  + \int_0^T \int_{\Omega} (\Lambda(\theta) + \varphi) \, {\bf u} \cdot \nabla_x \xi
  + \int_0^T \int_{\Omega} h(\theta) \Delta \xi\\
 \label{wf6}
  & \le - \int_0^T \int_{\Omega} \Bigg (\frac{\nu(\theta)}{\theta} |\nabla_x {\bf u}|^2
   + \frac{1}{\theta} |\nabla_x \mu|^2 + \frac{\kappa(\theta)}{\theta^2} |\nabla_x \theta|^2 \Bigg ) \xi
   - \int_{\Omega} (\Lambda(\theta_0) + \varphi_0)  \cdot \xi(0, \cdot),
\end{align}
holding for any $\xi \in \mathcal{C}^{\infty}_0([0,T) \times \Omega)$,
$\xi \ge 0$, and where we have set
\begin{equation}\label{defih}
  h(\theta)=\int_1^\theta
   \frac{\kappa(s)}s\,\dis
   = \log \theta + \frac{1}{\beta} (\theta^\beta - 1).
\end{equation}
\end{thm}
\begin{oss}\label{rego}
Let us notice that in case we could prove existence of a {\em sufficiently
smooth}\/ weak solution (in particular, regular enough
in order to integrate back by parts the terms in \eqref{wf5}),
then it would be possible to show that such a solution
also satisfies the ``standard'' form of the
heat equation \eqref{model5}. Equivalently, the entropy inequality
\eqref{entropy} would hold as an equality in that case. Hence, the current
notion of weak solution turns out to be compatible both
with Thermodynamics and also with the ``strong'' one.
\end{oss}


\section{A priori bounds}
\label{sec:bou}

The remainder of the paper is devoted to the proof of Theorem~\ref{thm:main}.
We start by briefly sketching our strategy.
In this section, we will prove some
{\it formal}\/ a-priori estimates holding for a hypothetical
quadruple $(\vu,\fhi,\mu,\theta)$
solving the ``strong'' formulation of the model stated in
Subsection~\ref{sec:pdes}. Actually, these estimates will mainly follow
as direct consequences of the Total energy balance~\eqref{energy0}
and Entropy inequality~\eqref{entropy}. Of course, to make
this procedure fully rigorous, one should rather work on a proper
regularization or approximation of the strong system
and prove that it admits at least one solution
being sufficiently smooth in order to comply with the estimates.
However, the system stated in Subsection~\ref{sec:pdes}
is rather complex and the related approximation, on
the one hand, would be particularly long and technical, and,
on the other hand, should not present particular difficulties, or
novelties, from the analytical viewpoint.
Indeed, the nonisothermal Navier-Stokes system given by
\eqref{momentumbal} and \eqref{theta} can be treated
along the lines developed in the monograph~\cite{FEI09},
while the regularization of the Cahn-Hilliard equation~\eqref{CH}
is completely standard since no singular or nonsmooth terms are involved.
For all these reasons, we decided to skip the details
of this argument and rather proceed formally.
We just quote the papers~\cite{A09, FFRS, FPR09, FRSZ, ls, ls1, lss} 
where more details of possible approximations of related
systems are given.

In Section \ref{WSS}, having the a-priori estimates at disposal,
we will then prove that any sequence
$({\bf u}_n, \varphi_n, \mu_n, \theta_n)$
complying with the bounds uniformly in $n$ admits
at least one limit point
$({\bf u}, \varphi, \mu, \theta)$ which solves the {\em weak}\/
formulation of the system (i.e., satisfies the conditions
stated in Definition~\ref{defweak}).
This procedure, which will be referred to as
{\it ``weak sequential stability''}\/ of families of solutions,
can be seen as a simplified version of the compactness argument
that one should use to remove some form of regularization
or approximation.


\subsection{Energy estimates}

Integrating the total energy balance \eqref{energy0},
we deduce the following a priori estimates:
\begin{eqnarray}
&& ||Q(\theta)||_{L^{\infty}(0,T;  L^1(\Omega))} \le c, \label{EN1}\\
&& ||{\bf u}||_{L^{\infty} (0,T; L^2(\Omega; \mathbb{R}^3))} \le c, \label{EN2}\\
&& ||F(\varphi)||_{L^{\infty}(0,T; L^1(\Omega))} \le c, \label{EN3}\\
&& {||\varphi||_{L^{\infty}(0,T; H^1(\Omega))} \le c.} \label{EN4}
\end{eqnarray}
By \eqref{kappacv}--\eqref{hpbetadelta}, \eqref{EN1} gives
\begin{equation}
\label{theta_cV}
||\theta||_{L^{\infty}(0,T; L^{\delta + 1}(\Omega))} \le c.
\end{equation}


{\subsection{Entropy estimates}}

Integrating the entropy inequality \eqref{entropy} both in space and
in time, and using \eqref{EN4}, \eqref{theta_cV}, and \eqref{incom},
we readily deduce
\begin{eqnarray}
&& {||\theta^{- 1/2} \nabla_x {\bf u}||_{L^2((0,T) \times \Omega; \mathbb{R}^{3 \times 3})} \le c,} \label{entropy1}\\
&&  {||\theta^{- 1/2} \nabla_x \mu||_{L^2((0,T) \times \Omega; \mathbb{R}^3)} \le c,} \label{entropy2}\\
&& \displaystyle {\int_0^T \int_{\Omega} \frac{\kappa(\theta)}{\theta^2} |\nabla_x \theta|^2 \le c.} \label{entropy3}
\end{eqnarray}
{By \eqref{kappacv}, \eqref{entropy3}, and due to the inequality
\[
  1 \le C \, \left (\frac{1}{x^2} + x^{\beta - 2} \right)
   \quad \forall x>0, 
\]
holding for $\beta \geq 2$ (cf.~\eqref{hpbetadelta}), we obtain on the one hand}
\begin{equation}
\label{regtheta1}
||\nabla_x \theta||_{L^2((0,T) \times \Omega; \mathbb{R}^3)} \le c.
\end{equation}
On the other hand, again from \eqref{kappacv} and \eqref{entropy3} we infer
\begin{align}
 &||\nabla_x (\log(\theta))||_{L^2((0,T) \times \Omega; \mathbb{R}^3)} \le c, \label{regtheta2}\\
 &||\nabla_x (\theta^{\frac{\beta}{2}})||_{L^2((0,T) \times \Omega; \mathbb{R}^3)} \le c. \label{regtheta3}
\end{align}
{At this point, \eqref{theta_cV} and \eqref{regtheta3} together with a generalized
version of Poincar\'e's inequality yield}
\[
||\theta^{\beta/2}||_{L^2(0,T; L^6(\Omega))} \le c,
\]
{from which we deduce}
\begin{equation}
\label{reg_theta}
||\theta||_{L^{\beta}(0,T; L^{3 \beta}(\Omega))} \le c.
\end{equation}


\subsection{Temperature estimates}

We now integrate the temperature equation~\eqref{theta}
in space and time, with the aim of getting a $L^2$-bound for the two
quadratic terms on the right hand side. Note that, at this level
(or, to be more precise, in the approximation we decided to skip),
it is crucial to have the strong relation~\eqref{theta} (or an
approximated version of it ``containing'' the same information)
at disposal. In other words, this procedure
cannot be reproduced by workly directly
on weak solutions since the total energy
balance~\eqref{wf5} alone is not sufficient.

That said, we need to control the terms on the left hand side
of \eqref{theta}. First of all, \eqref{EN1} brings
\[
\int_{\Omega} \int_0^T (c_V(\theta) \theta_t)  \le c.
\]
Moreover, by boundary conditions and incompressibility,
\[
  \int_0^T \int_{\Omega} \dive (\kappa(\theta) \nabla_x \theta)
   = \int_0^T \int_{\Omega} {\bf u} \cdot \nabla_x Q(\theta) = 0.
\]
Finally, using the first equation of \eqref{CH}, we are able to deduce
\[
 \int_0^T \int_{\Omega} \theta \frac{D \varphi}{D t}
  = \int_0^T \int_{\Omega} \theta (\varphi_t + {\bf u} \cdot \nabla_x \varphi)
  = - \int_0^T \int_{\Omega} \nabla_x \mu \cdot \nabla_x \theta
  \le \frac12 \int_0^T \int_{\Omega} | \nabla_x \mu |^2
   + \frac12 \int_0^T \int_{\Omega} | \nabla_x \theta |^2
\]
and we need to control the last two terms. Actually,
the first one is absorbed by the last term on the
\rhs\ of \eqref{theta}, while the second one is controlled
thanks to \eqref{regtheta1}.

Hence, from the right hand side of \eqref{theta} we can ``read''
the a priori estimates
\begin{eqnarray}
&& ||D {\bf u}||_{L^2((0,T) \times \Omega; \mathbb{R}^{3\times 3})} \le c, \label{Duno}\\
&& ||\nabla_x \mu||_{L^2((0,T) \times \Omega; \mathbb{R}^3)} \le c. \label{Ddue}
\end{eqnarray}


{\subsection{Consequences}}

\subsubsection{Higher regularity for $\varphi$ and $F'(\varphi)$}

First of all, thanks to \eqref{EN4} and the classical Sobolev embedding theorems, we deduce
\begin{equation}
\label{varphiinfty6}
||\varphi||_{L^{\infty}(0,T; L^6(\Omega))} \le c,
\end{equation}
and therefore
\begin{equation}
\label{F-infty2}
||F'(\varphi)||_{L^{\infty}(0,T; L^2(\Omega))} \le c.
\end{equation}
Let us note that here it is essential to have a cubic growth in 3D for the potential $F'$.
Now we choose $\xi = \varphi$ in \eqref{wf3}, we test \eqref{wf4} by
$\Delta \varphi$, and take the sum. Some terms cancel out
due to \eqref{incom} and our choice of boundary conditions.
Hence, using in particular \eqref{theta_cV}, \eqref{varphiinfty6} and \eqref{F-infty2}
for treating the other terms, it is not difficult to get the
additional regularity on $\varphi$:
\begin{equation}
\label{phi1}
  ||\Delta \varphi||_{L^2((0,T) \times \Omega)} \le c \,\,
     \Rightarrow  ||\varphi||_{L^2(0,T; H^2(\Omega))} \le c,
\end{equation}
where the classical regularity theorems for elliptic equations have
also been used.
\\
Now, we would like to show that
\begin{equation}
\label{nablaF22}
||\nabla_x (F'(\varphi))||_{L^2((0,T) \times \Omega; \mathbb{R}^3)} \le c.
\end{equation}
To this aim, we use classical interpolation inequalities. Indeed,
\begin{eqnarray*}
\varphi \in L^{\infty}(0,T; H^1(\Omega)) \cap L^2(0,T; H^2(\Omega))
 &\hookrightarrow & L^{2/\vartheta} (0,T; H^{1+ \vartheta}(\Omega)) \\
&\hookrightarrow & L^{2/\vartheta} (0,T; L^{\frac{6}{1 - 2 \vartheta}}(\Omega)),
\end{eqnarray*}
provided that $\vartheta \in (0, 1/2)$, whence
\begin{equation}
\label{uno}
F''(\varphi) = 3 \varphi^2 - 1 \in L^{1/\vartheta} (0,T; L^{\frac{3}{1 - 2 \vartheta}}(\Omega)).
\end{equation}
On the other hand
\begin{eqnarray}
\nabla_x \varphi \in L^{\infty}(0,T; L^2(\Omega; \mathbb{R}^3)) \cap L^2(0,T; H^1(\Omega; \mathbb{R}^3)) &\hookrightarrow& L^{2/\vartheta} (0,T; H^{\vartheta}(\Omega; \mathbb{R}^3))  \label{due}\\
&\hookrightarrow & L^{2/\vartheta} (0,T; L^{\frac{6}{3 - 2 \vartheta}}(\Omega; \mathbb{R}^3)). \nonumber
\end{eqnarray}
At this point, combining \eqref{uno} and \eqref{due} we need to find $\vartheta < 1/2$ such that
\[
\vartheta + \frac{\vartheta}{2} = \frac{1}{2} \,\,\, \textnormal{and} \,\,\, \frac{1- 2 \vartheta}{3} + \frac{3 - 2 \vartheta}{6} = \frac{1}{2}
\]
and this leads to $\vartheta = 1/3$. For this value of $\vartheta$ we then have
\[
F''(\varphi) \in L^3 (0,T; L^9(\Omega)), \qquad \nabla_x \varphi \in L^6(0,T; L^{18/7}(\Omega; \mathbb{R}^3)),
\]
whence \eqref{nablaF22}.


\subsubsection{Further regularity for $\varphi$, $\nabla_x \varphi$ and $\mu$}

First of all, integrating in space the second equation of \eqref{CH}, using \eqref{theta_cV}, \eqref{BC3} and \eqref{F-infty2}, we arrive at
\[
||\overline{\mu}||_{L^\infty(0,T)} \le c,
\]
which implies, due to \eqref{Ddue} and the Poincar\'e-Wirtinger inequality,
\begin{equation}
\label{mu}
||\mu||_{L^2(0,T; H^1(\Omega))} \le c.
\end{equation}
At this point, in view of {\eqref{regtheta1}, \eqref{nablaF22} and \eqref{mu}}, we can interpret the second equation of \eqref{CH} as $-\Delta \varphi = f:=\mu+\theta-F'(\varphi)$,
where, by direct comparison, $f \in L^2(0,T; H^1(\Omega))$. By
the classical regularity theorems for elliptic equations this implies
\[
  ||\varphi||_{L^2(0,T; H^3(\Omega))} \le c.
\]
Combining this relation with \eqref{EN4}, we have
\begin{equation}\label{Enablafi}
  ||\nabla_x \varphi||_{L^2(0,T; H^2(\Omega;\RR^3)) \cap L^{\infty}(0,T; L^2(\Omega;\RR^3))} \le c.
\end{equation}
Using interpolation and the Gagliardo-Nirenberg inequalities
(cf., e.g., \cite[p.~125]{nier})
we then deduce
%
%
\begin{equation}
\label{laplaphi}
  ||\varphi||_{L^{\frac{14}{5}}(0,T; W^{2,{\frac{14}{5}}}(\Omega))},
   \quad \|\nabla\varphi\|_{L^{\frac{14}{3}}((0,T)\times\Omega;\RR^3)}\leq c,
\end{equation}
%
as well as
\begin{equation}
\label{nablafi}
||\nabla \varphi||_{L^3(0,T; W^{1,{\frac{18}{7}}}(\Omega; \mathbb{R}^3))} \le c
  \Rightarrow ||\nabla \varphi||_{L^3(0,T; L^{18}(\Omega; \mathbb{R}^3))} \le c.
\end{equation}
To conclude, it remains to show that
\begin{equation}
\label{reg_mu}
||\mu||_{L^{{\frac{14}{5}}}((0,T) \times \Omega)}\leq c.
\end{equation}
We proceed by comparing terms in the second of \eqref{CH}.
Indeed, due to \eqref{F-infty2} and \eqref{nablaF22},
\[
F'(\varphi) \in L^{\infty}(0,T; L^2(\Omega)) \cap L^2(0,T; H^1(\Omega)) \hookrightarrow L^{\frac{10}{3}}((0,T) \times \Omega) \hookrightarrow L^{\frac{14}{5}}((0,T) \times \Omega).
\]
On the other hand,
\begin{equation}\label{XY}
\theta \in L^{\infty}(0,T; L^{1 + \delta}(\Omega)) \cap L^{\beta}(0,T; L^{3 \beta}(\Omega)) \hookrightarrow L^{p_{\beta, \delta}}((0,T) \times \Omega) \hookrightarrow L^{\frac{14}{5}}((0,T) \times \Omega)
\end{equation}
as long as
\begin{equation}\label{XYZ}
p_{\beta,\delta} > 3 > \frac{14}{5},
\end{equation}
which is true thanks to \eqref{hpbetadelta} (cf.~also \eqref{cost_A}).
Recalling also \eqref{laplaphi}, we then get \eqref{reg_mu}.


\section{Weak sequential stability}
\label{WSS}

In this section, we assume to have a sequence
$({\bf u}_n, \varphi_n, \mu_n, \theta_n)$ of solutions satisfying
(a proper approximation of) the strong system of Subsection~\ref{sec:pdes}.
Then, by virtue of the argument developed in the previous part,
we can assume that this family complies with the proved
a priori bounds uniformly with respect to $n$. Our aim is
showing, by weak compactness arguments, that at least a
subsequence converges in a suitable way to a weak solution
to our problem (i.e.~to a limit quadruple $({\bf u}, \varphi, \mu, \theta)$
satisfying the statement given in Definition~\ref{defweak}).
Actually, to further simplify the notation, we intend that
all the convergence relations appearing in the following
are to be considered up to the extraction of (not relabelled)
subsequences.

That said, collecting the bounds proved before, we have
\begin{eqnarray}
&& {\bf u}_n \rightarrow {\bf u} \,\,\, \textnormal{weakly star in $L^{\infty}(0,T; L^2(\Omega; \mathbb{R}^3)) \cap L^2(0,T; H^1(\Omega; \mathbb{R}^3))$},  \label{limit1}\\
&& {\varphi}_n \rightarrow {\varphi} \,\,\, \textnormal{weakly star in $L^{\infty}(0,T; H^1(\Omega)) \cap L^2(0,T; H^3(\Omega))$}, \label{limit2}\\
&& \mu_n \rightarrow \mu \,\,\, \textnormal{weakly in $L^2(0,T; H^1(\Omega)) \cap L^{\frac{14}{5}}((0,T) \times \Omega)$}, \label{limit3}\\
&& \theta_n \rightarrow \theta \,\,\, \textnormal{weakly star in $L^{\infty}(0,T; L^{\delta+1}(\Omega)) \cap L^{\beta}(0,T; L^{3 \beta}(\Omega))$}, \label{limit4}\\
&& \theta_n \rightarrow \theta \,\,\, \textnormal{weakly in $L^2(0,T; H^1(\Omega))$}, \label{limit4b}
\end{eqnarray}
where we recall that $\beta$ and $\delta$ fufill assumption \eqref{hpbetadelta}.

Next, we need an estimate for the pressure $p_n$.
%
%
To achieve it, we will follow the approach devised in \cite{BFM09}
for dealing with the Navier-Stokes-Fourier system.
Referring to that paper for more details, a formal
way to get integrability of $p_n$ consists in computing
the divergence of \eqref{momentum}, which gives rise to
\begin{equation} \label{ellprob}
  \Delta p_n = \dive \dive \big( \mathbb{S}_n - {\bf u}_n \otimes {\bf u}_n
    - \Grad \fhi_n \otimes \Grad \fhi_n \big).
\end{equation}
In other words, $p_n$ solves, at least formally, some kind of elliptic problem.
Hence, an estimate for it can be proved by relying on
suitable regularity theorems.
To be more precise, \eqref{ellprob} has to be interpreted in a ``very weak''
sense. Namely, $p_n$ turns out to satisfy the integral identity
\begin{equation} \label{press}
  \int_\Omega p_n \Delta \xi
   =  \int_{\Omega} \big( \mathbb{S}_n - {\bf u}_n \otimes {\bf u}_n
         - \Grad \fhi_n \otimes \Grad \fhi_n \big) :  \Grad \Grad \xi,
\end{equation}
for any test function $\xi\in H^2(\Omega)$ with
$\nabla_x\xi \cdot {\bf n}|_{\Gamma} = 0$.

The above formulation incorporates the boundary conditions,
though in a way which is not completely obvious.
In particular, one uses in an
essential way the complete slip conditions~\eqref{BC1}
for $\vu$ and the no-flux conditions for the other variables.
Indeed, the Neumann condition \eqref{BC3} allows us to deal
with the extra stress $- \Grad \fhi_n \otimes \Grad \fhi_n$. To
see this at least formally, one multiplies \eqref{momentumbal}
by $\Grad \xi$. Then, thanks to \eqref{BC3} and \eqref{BC1}, it is possible to
integrate by parts and get \eqref{press} without the occurrence of additional
boundary terms. This procedure can be made rigorous by following
closely the lines of the approximation argument
described in~\cite[Sec.~4]{BFM09}, to which
we refer the reader for more details.

That said, we want to apply suitable elliptic regularity theorems
to \eqref{ellprob} (or, more precisely, to its weak formulation
\eqref{press}). To this aim, we first establish some bounds
on the terms appearing on the right hand side.

Firstly, we notice that
\begin{equation}
\label{embed_u}
L^{\infty}(0,T; L^2(\Omega; \mathbb{R}^3)) \cap L^2(0,T; H^1(\Omega; \mathbb{R}^3))
\hookrightarrow L^{10/3}((0,T) \times \Omega; \mathbb{R}^3),
\end{equation}
continuously. Hence, \eqref{limit1} entails
\begin{equation}
\label{bound_u}
  \big\| |{\bf u}_n|^3  \big\|_{L^{1^+}((0,T) \times \Omega)} \le c,
\end{equation}
as well as
\begin{equation}
\label{aaa}
||{\bf u}_n \otimes {\bf u}_n||_{L^{{\frac{5}{3}}}((0,T)\times \Omega; \mathbb{R}^{3 \times 3})} \le c.
\end{equation}
Here and in the sequel the notation $1^+$ indicates a proper
exponent strictly greater than 1. Due to \eqref{Duno}
and the last \eqref{kappacv}, we also have
\begin{equation}
\label{bbb}
||\mathbb{S}_n||_{L^{{2}}((0,T)\times \Omega; \mathbb{R}^{3 \times 3})} \le c.
\end{equation}
On the other hand, thanks to \eqref{laplaphi}, 
we have
\begin{equation}\label{ccc}
  \big\| \Grad \fhi_n \otimes \Grad \fhi_n \big\|_{ L^{\frac{7}{3}}
      ((0,T)\times \Omega; \mathbb{R}^{3\times 3} )} \le c.
\end{equation}
Then, collecting \eqref{bound_u}--\eqref{ccc} and using elliptic regularity
in \eqref{press} (again, the precise details are given in~\cite{BFM09}),
we deduce
\begin{equation}
\label{pressione}
  p_n \rightarrow p \qquad \textnormal{weakly in {$L^{\frac{5}{3}}((0,T)\times \Omega)$}}.
\end{equation}
Consequently, a comparison of terms in \eqref{momentum} gives
\[
||({\bf u}_{n})_t||_{L^{\frac{5}{3}}(0,T; X)} \le c,
\]
where $X$ is a suitable Sobolev space of negative order.
By the Aubin-Lions lemma, \eqref{limit1} and \eqref{bound_u}, we then get
\begin{equation}
\label{strongu}
{\bf u}_n \rightarrow {\bf u} \qquad \textnormal{strongly in $L^{3^+}((0,T) \times \Omega; \mathbb{R}^3$)}.
\end{equation}
Next, we observe that ${\bf u}_n \cdot \Grad \varphi_n$
is bounded in $ L^2((0,T)\times\Omega)$, due to \eqref{EN2} and \eqref{Enablafi}. 
Hence, recalling \eqref{limit3} and comparing terms in the first \eqref{CH} we infer
\[
  ||(\varphi_n)_t||_{L^2(0,T; (H^1(\Omega))')} \le c.
\]
Combining this with \eqref{limit2} and applying once more the Aubin-Lions lemma, we get
\begin{equation}
\label{strongphi}
\varphi_n \rightarrow \varphi \qquad \textnormal{strongly in $L^2(0,T; H^{3-\sigma}(\Omega))
   \cap \mathcal{C}^0([0,T]; L^2(\Omega))$},
\end{equation}
for all $\sigma > 0$.

Now, let us assume that, for every $n\in\mathbb{N}$,   
the approximate solution $({\bf u}_n, \varphi_n, \mu_n, \theta_n)$ 
fulfills the strong system stated in Subsec.~\ref{subs:strong}
(actually, its hypothetical approximation).
Then, multiplying by suitable
test functions and integrating by parts, it readily
follows that $({\bf u}_n, \varphi_n, \mu_n, \theta_n)$
also complies with Definition~\ref{defweak}. In particular, 
starting from the ``strong'' formulation \eqref{energy0}
of the total energy balance (assumed to hold at the $n$-level),
accounting for \eqref{towf11} and performing suitable
integrations by parts, it is not difficult to deduce
its weak counterpart
\begin{align} \nonumber
  & \partial_t \left( \frac{1}{2}|\bu_n|^2 + e_n \right)
   + \dive \left( \frac{1}{2} \bu_n |\bu_n|^2 + e_n \bu_n \right)+\dive(p_n{\bf u}_n)
   - \Delta \hat{\kappa}(\theta_n)
   - \dive ( \mathbb{S}_n \bu_n )\\
 \label{energy1}
  & \mbox{}~~~~~
   - \Delta \mu_n^2
   + \dive \big( (\vu_n \cdot \Grad \fhi_n ) \Grad \fhi_n \big)
   + \dive ( \Grad \Grad \fhi_n \Grad \mu_n )
   - \dive \dive ( \Grad \mu_n \otimes \Grad \fhi_n )
    = 0,
\end{align}
which holds at least in the sense of distributions,
and where $e_n$ is defined as 
\begin{equation}\label{defen}
  e_n := F(\varphi_n) + \frac{|\nabla_x \varphi_n|^2}{2} + Q(\theta_n). 
\end{equation}
Let us note that \eqref{energy1} reduces to \eqref{wf5} after testing it
by $\xi \in C_0^\infty([0,T)\times\Omega)$.
It is worth remarking that, in the approximation,
$({\bf u}_n, \varphi_n, \mu_n, \theta_n)$ also
satisfies the strong form~\eqref{entropy} of the entropy 
relation.

Then, we can see what happens as we let $n\nearrow \infty$.
Actually, on account of the properties proved above,
it is a standard matter to see that both the Cahn-Hilliard system
\eqref{wf3}--\eqref{wf4} and the momentum equation \eqref{wf1}
pass to the desired limits. Hence, to conclude the proof,
we have to prove \eqref{wf5} and \eqref{wf6}
in the limit. This issue is a bit more delicate and is 
dealt with in the next two subsections.


\subsection{Limit of the total energy balance}
\label{lim:bal}

Here, we shall take the limit $n\searrow \infty$ in~\eqref{energy1}. 
Let us start considering the term $e_n{\bf u}_n$.
We can first notice that, by \eqref{defQ} and \eqref{XY}, there follows 
\[
  \|Q(\theta_n)\|_{L^{\frac{p_{\beta,\delta}}{\delta+1}}((0,T)\times\Omega)}\leq c,
\]
where $p_{\beta, \delta}/(\delta+1)>5/3$ due to \eqref{cost_A} and \eqref{hpbetadelta}.  
Hence, using \eqref{bound_u}, we obtain 
\[
  \|Q(\theta_n) \, {\bf u}_n\|_{L^{1^+}((0,T) \times \Omega; \mathbb{R}^3)} \le c.
\]
Thus, thanks also to \eqref{F-infty2} and \eqref{nablafi}, we get 
\[
  \|e_n{\bf u}_n\|_{L^{1+}((0,T)\times \Omega; \RR^3)}\leq c. 
\]
Moreover, due to  \eqref{bound_u} and \eqref{pressione}, we obtain further that
\[
  || p_n \, {\bf u}_n||_{L^{1^+}((0,T) \times \Omega; \mathbb{R}^3)} \le c.
\]
Next, let us turn our attention to the term $\hat{\kappa}(\theta_n)$.
In order to take its limit we need (cf.~\eqref{kciapo})
the uniform bound
\begin{equation}
\label{Eteta}
||\theta_n||_{L^{(\beta + 1)^+}((0,T) \times \Omega)}\le c
\end{equation}
together with pointwise (a.e.) convergence, which will be shown
below. Thanks to \eqref{XY}, \eqref{Eteta} holds whenever
\begin{equation}
\label{gamma_hp}
{p_{\beta, \delta}} > \beta + 1 \,\, \Leftrightarrow \,\, \delta > \frac{1}{2}.
\end{equation}
Recalling that thermodynamical consistency of the model requires $\delta < 1$,
assumption \eqref{hpbetadelta} is then fully justified.

Finally, let us deal with the remaining terms in the total
energy balance \eqref{wf5}. 
Recalling that $\mathbb{S}_n = \nu(\theta_n) D {\bf u}_n$,
and using \eqref{Duno} and \eqref{bound_u}, we also obtain
\[
  ||\mathbb{S}_n {\bf u}_n||_{L^{1^+}((0,T) \times \Omega; \mathbb{R}^3)} \le c.
\]
Now, by \eqref{reg_mu} we have
\begin{equation}\label{esnew1}
  ||\mu_n^2||_{L^{1^+}((0,T) \times \Omega)}\le c.
\end{equation}
Next, thanks to \eqref{bound_u} and \eqref{laplaphi}, we infer
\begin{equation}\label{esnew2}
  ||({\bf u}_n \cdot \nabla_x \varphi_n) \nabla_x \varphi_n ||%
   _{L^{1^+}((0,T) \times \Omega;\mathbb{R}^3)}\le c.
\end{equation}
Moreover, \eqref{mu} and \eqref{laplaphi} lead to
\begin{equation}\label{esnew3}
 \big\| ( \nabla_x \nabla_x \varphi_n ) \nabla_x \mu_n  \big\|_{L^{1^+}((0,T)\times \Omega; \mathbb{R}^3)} +\big\| \nabla_x \mu_n \otimes \nabla_x \varphi_n\big\|_{L^{1^+}((0,T) \times \Omega; \mathbb{R}^{3 \times 3})}\le c.
\end{equation}
%
%
%
%
%
%
%
At this point, collecting the previous estimates and comparing terms in \eqref{energy1},
we obtain
\[
  \left | \left | \partial_t \left ( \frac{1}{2} |{\bf u}_n|^2 + e_n\right )\right | \right |_{L^{1^+}(0,T; X)}
   \le c,
\]
where $X$ is, again, some Sobolev space of negative order. On the other hand,
a direct computation based on  the estimates \eqref{Eteta}, \eqref{regtheta1},
\eqref{bound_u}, \eqref{laplaphi}, and \eqref{nablafi} yields
\[
  \left | \left |\nabla_x \left ( \frac{1}{2} |{\bf u}_n|^2
   + e_n\right )\right | \right |_{L^{1^+}((0,T) \times \Omega; \mathbb{R}^3)} \le c.
\]
Hence, we can use once more the Aubin-Lions lemma to conclude that
\[
  \left (\frac{1}{2} |{\bf u}_n|^2 + e_n \right )
    \rightarrow \left (\frac{1}{2} |{\bf u}|^2 + e \right ) \qquad
   \textnormal{strongly in ${L^{1^+}((0,T) \times \Omega)}$}.
\]
Actually, the limit of $|{\bf u}_n|^2$ is identified as
$|{\bf u}|^2$ thanks to the strong convergence \eqref{strongu},
while the limit $e$ of $e_n$ (cf.~\eqref{defen}) 
still needs to be identified in terms of $\theta$ and $\varphi$.
To this aim, we need strong convergence of $\theta_n$ in $L^p((0,T) \times \Omega)$
for some $p$ and, to achieve it, we use a monotonicity argument like in the paper \cite{FRSZ}.

Actually, from \eqref{strongphi} we infer
\[
  F(\varphi_n) + \frac{|\nabla_x \varphi_n|^2}{2}
   \rightarrow F(\varphi) + \frac{|\nabla_x \varphi|^2}{2},
   \quad\hbox{strongly, say, in } L^{1^+}((0,T)\times\Omega).
\]
Hence, taking a couple of indexes $n$ and $m$ and using the previous
information with the strong $L^1$-convergence of $e_n \to e$,
we get
\begin{align*}
  & \int_0^T \int_{\Omega} (Q(\theta_n) - Q(\theta_m),  \textnormal{sign} (\theta_n - \theta_m))
    = \int_0^T \int_{\Omega} (e_n - e_m, \textnormal{sign} (\theta_n - \theta_m)) \\
  & \mbox{}~~~~~
    - \int_0^T \int_{\Omega} \left (F(\varphi_n) + \frac{|\nabla_x \varphi_n|^2}{2}
      - F(\varphi_m) - \frac{|\nabla_x \varphi_m|^2}{2},
    \textnormal{sign} (\theta_n - \theta_m) \right )\rightarrow 0.
\end{align*}
Due to monotonicity of $Q$, we then get that $Q(\theta_n)$ is a Cauchy sequence
in $L^1((0,T)\times\Omega)$. Hence it converges strongly and,
up to a subsequence, almost everywhere, to some limit $q$. Combining this
fact with \eqref{limit4} and using a generalized form of Lebesgue's theorem, we
get that $\theta=\big((\delta+1)q\big)^{1/{(\delta+1)}}$
(cf.~\eqref{defQ}). Moreover, we obtain
\begin{equation}
\label{conv_theta}
  \theta_n \rightarrow \theta \qquad \textnormal{strongly in $L^p((0,T) \times \Omega)$
  for all $p \in [1, {p_{\beta, \delta}})$},
\end{equation}
where $p_{\beta, \delta}$ has been introduced in \eqref{cost_A}. This in particular implies that
\begin{eqnarray*}
&& \kappa(\theta_n) \rightarrow \kappa(\theta) \qquad \textnormal{strongly in $L^p((0,T) \times \Omega)$ for all $p \in \bigg[1, \frac{{p_{\beta, \delta}}}{\beta}\bigg)$},\\
&& c_V(\theta_n) \rightarrow c_V(\theta) \qquad \textnormal{strongly in $L^p((0,T) \times \Omega)$ for all $p \in \bigg[1, \frac{{p_{\beta, \delta}}}{\delta}\bigg)$}.
\end{eqnarray*}
In view of the above discussion,
all terms in the first row
of \eqref{energy1} pass to the desired limits.
On the other hand, the last three terms in the second row
can be managed thanks to (\ref{esnew1}--\ref{esnew3}). 
Finally, combining \eqref{strongphi} with \eqref{conv_theta}
and comparing terms in \eqref{wf4}, we obtain that
\begin{equation}\label{conv_mu}
  \mu_n \rightarrow \mu = - \Delta \fhi + F'(\fhi) - \theta
  \qquad \textnormal{strongly in $L^2((0,T) \times \Omega)$}.
\end{equation}
Hence, $\Delta \mu_n^2 \to \Delta \mu^2$ at least in the sense
of distributions. This allows us to take the limit of
\eqref{energy1}, which immediately reduces to
\eqref{wf5} after testing by $\xi\in \mathcal{C}_0^\infty([0,T)\times\Omega)$
and integrating by parts.


\subsection{Proof of the entropy inequality}
\label{lim:entropy}

To conclude the proof of Theorem~\ref{thm:main}, we need
to prove the entropy production inequality~\eqref{wf6}.
As noted above, we can assume that a stronger relation,
i.e.~\eqref{entropy}, holds at the $n$-level
and we aim at taking its (supremum) limit as $n\to \infty$. 
In particular, this will give rise to the $\geq$ sign in \eqref{wf6}.

For the reader's convenience, we start by reporting
the statement of a useful lower semicontinuity result
due to {\sc A.D.~Ioffe} \cite{Ioffe}:
\begin{thm}
\label{ioffe}
Let $\mathcal{O}\subset \mathbb{R}^d$ a smooth bounded open set and
$f: \mathcal{O} \times \mathbb{R}^n \times \mathbb{R}^m \rightarrow [0, + \infty]$,
$d ,n, m \ge 1$,
be a measurable non-negative function such that
\begin{eqnarray}
&& f(x, \cdot, \cdot) \,\,\, \textnormal{is lower semicontinuous on $\mathbb{R}^n \times \mathbb{R}^m$ for every $x \in \mathcal{O}$}, \label{ioffe1}\\
&& f(x, u, \cdot) \,\,\, \textnormal{is convex on $\mathbb{R}^m$ for every $(x,u) \in \mathcal{O} \times \mathbb{R}^n$.} \label{ioffe2}
\end{eqnarray}
Let also $(u_k, v_k), (u,v): \mathcal{O} \rightarrow \mathbb{R}^n \times \mathbb{R}^m$
be measurable functions such that
\[
u_k(x) \rightarrow u(x) \,\,\, \textnormal{in measure in $\mathcal{O}$}, \qquad v_k \rightharpoonup v \,\,\, \textnormal{weakly in $L^1(\mathcal{O}; \mathbb{R}^m)$}.
\]
Then,
\[
\liminf_{k \rightarrow + \infty} \int_{\mathcal{O}} f(x, u_k(x), v_k(x)) \ge \int_{\mathcal{O}} f(x, u(x), v(x)).
\]
\end{thm}
\noindent%
To begin, we test \eqref{entropy} at the $n$-level by a nonnegative
test function $\xi \in C_0^\infty([0,T)\times \Omega)$, as specified
in the statement of Theorem~\ref{thm:main},
and we integrate by parts. What we get is exactly the $n$-version of
\eqref{wf6} (with the equal sign). In particular, integration by
parts gives rise to the function $h(\theta_n)$ (cf.~\eqref{defih}).
Thanks to \eqref{regtheta2} and \eqref{conv_theta}
(cf.~also Remark~\ref{thetapos} below), it is not
difficult to see that
\begin{equation}\label{conv:h}
  h(\theta_n) \rightarrow h(\theta),
  \qquad \textnormal{say, strongly in $L^{1^+}((0,T) \times \Omega)$},
\end{equation}
whence in particular
\[
  \int_0^T \int_{\Omega} h(\theta_n) \Delta \xi
   \to \int_0^T \int_{\Omega} h(\theta) \Delta \xi.
\]
To deal with the other terms on the left hand side
of \eqref{wf6}$_n$, we observe that, due to \eqref{conv_theta}
and \eqref{Lam},
\[
  \Lambda(\theta_n) \rightarrow \Lambda(\theta)
   \qquad \textnormal{strongly in $L^p((0,T) \times \Omega)$
   for all $p \in \bigg [1, \frac{p_{\beta,\delta}}{\delta}\bigg )$},
\]
whence, recalling \eqref{strongu},
\[
  \Lambda(\theta_n) \vu_n \rightarrow \Lambda(\theta) \vu,
   \qquad \textnormal{say, strongly in $L^1((0,T) \times \Omega;\mathbb{R}^3)$}.
\]
The remaining terms on the \lhs\ of \eqref{wf6} are simpler to treat.
To deal with the \rhs, we recall
\eqref{limit1}, \eqref{limit3}, and \eqref{limit4b},
and use Ioffe's theorem. This gives
\begin{eqnarray}
  && \int_0^T \int_{\Omega} \left ( \frac{\xi}{\theta}
   \left (\nu(\theta) |\nabla_x {\bf u}|^2 + |\nabla_x \mu|^2
     + \frac{\kappa(\theta)}{\theta} |\nabla_x\theta|^2 \right ) \right ) \nonumber\\
 \label{fatou}
  & \le& \liminf_{n \rightarrow \infty} \int_0^T \int_{\Omega} \left ( \frac{\xi}{\theta_n}
    \left (\nu(\theta_n) |\nabla_x {\bf u}_n|^2 + |\nabla_x \mu_n|^2
       + \frac{\kappa(\theta_n)}{\theta_n} |\nabla_x\theta_n|^2 \right ) \right ).
\end{eqnarray}
Hence, \eqref{wf6} passes to the {\it supremum}\/ limit
$n\to \infty$. This concludes the proof of Theorem~\ref{thm:main}.
%
%
%
\begin{oss}\label{thetapos}
 In order for our estimates to make sense, we implicitly assumed 
 in the course of the proof the temperature
 $\vt$ to be (almost everywhere) positive. This fact is used
 in several estimates which, otherwise, would not make
 sense. Positivity of $\vt_n$ should be shown, indeed,
 at the $n$-level, i.e., for the hypothetical
 regularized problem which we decided not to detail here.
 Actually, for regularized solutions, which are usually smoother,
 it is often possible to prove some stronger property
 (like $\theta_n(x) \ge c_n >0$ a.e.),
 by applying a suitable maximum principle.
 We cannot give here a proof of this fact, since this
 would require to provide the details of the regularization. However,
 we can at least show that, if $\vt_n$ is almost
 everywhere positive, and satisfies the estimates given
 in Section~\ref{sec:bou}, then positivity is preserved in
 the limit. To see this, we first notice that, by \eqref{regtheta2},
 \begin{equation}\label{10}
   \| \Grad\log\vt_n \|_{L^2((0,T)\times \Omega;\mathbb{R}^3)} \le c.
 \end{equation}
 Then, using that the mean value of $\varphi_n$ is conserved, i.e.,
 \[
   \overline{\varphi_n}(t)=\overline{\varphi_n}(0)\quad\hbox{ for every } t\in (0,T)
 \]
 and integrating \eqref{entropy} both in space and in time,
 we readily obtain
 $$
   \io \Lambda(\vt_n(t)) \ge \io \Lambda(\vt_0)\quad \hbox{for a.e. }t\in (0,T)
 $$
 or, equivalently,
 \begin{equation}\label{11}
   \io \vt_n^\delta(t) \ge \io \vt_0^\delta\quad \hbox{for a.e. }t\in (0,T).
 \end{equation}
 Combining \eqref{10} and \eqref{11}, it is not difficult to deduce
 \begin{equation}\label{11.2}
   \| \log\vt_n \|_{L^2(0,T;H^1(\Omega))} \le c.
 \end{equation}
 This fact, together with the pointwise (a.e.)~convergence of $\vt_n$,
 implies
 \begin{equation}\label{15}
  \log \vt_n \to \log \vt,
   \ \ \text{say, strongly in }L^{2^-}((0,T)\times\Omega).
 \end{equation}
 In particular, $\vt > 0$ almost everywhere also in the limit.
\end{oss}


%


\bibliographystyle{amsplain}

\end{document}